\documentclass[acmlarge]{acmart}

\usepackage{booktabs} % For formal tables

\usepackage[ruled]{algorithm2e} % For algorithms

\SetAlFnt{\small}
\SetAlCapFnt{\small}
\SetAlCapNameFnt{\small}
\SetAlCapHSkip{0pt}
\IncMargin{-\parindent}

\begin{document}
% pierwsza strona pusta
% ~\\
% \newpage
% Title portion
\title{Functionals and hardware} 
\author{Stanislaw Ambroszkiewicz }
\orcid{}
\affiliation{%
  \institution{Institute of Computer Science, Polish Academy of Sciences, Jana Kazimierza 5, PL-01-248 Warsaw, Poland }
  \streetaddress{Jana Kazimierza 5}
  \city{Warsaw } 
  \state{}
  \postcode{PL-01-248}
  \country{Poland}}
%\affiliation{%
 % \institution{Siedlce University of Natural Sciences and Humanities}
 % \city{Siedlce}
 % \country{Poland}}

\begin{abstract}
Functionals are an important research subject in Mathematics and Computer Science as well as a challenge in Information Technologies where the current programming paradigm states that only symbolic computations are possible on higher order objects, i.e. functionals are terms, and computation is term rewriting. The idea explored in the paper is that functionals correspond to generic mechanisms for management of connections in  arrays  consisting of first order functional units.  Functionals are higher order abstractions that are useful for the management of such large arrays. Computations on higher order objects comprise dynamic configuration of connections between first order elementary functions in the arrays. Once the functionals are considered as the generic mechanisms, they have a grounding in hardware.  
A conceptual framework for constructing such mechanisms is presented, and their hardware realization is discussed.  
\end{abstract}


% We no longer use \terms command
%\terms{ }


%\thanks{This work is supported by the National Science Foundation, under grant CNS-0435060, grant CCR-0325197 and grant EN-CS-0329609.}

\maketitle

%\renewcommand{\shortauthors}{G. Zhou et al.}

\section{Functionals in mathematical logic } 

Functionals are higher order functions that take functions (as well as functionals) as arguments. The prominent examples of functionals in Mathematics are: limit of a sequence of numbers, derivative and integration. Composition of two function may be abstracted to a functional. Application of a function to its argument may be considered as a functional if abstracted from the function and the argument. Primitive recursion schema, where function is a parameter, may be viewed as a functional if the parameter is considered as argument. Fix-point operator is a functional. Generally, functionals are higher order abstractions that are not easy to comprehend. Since they are crucial in computing theory, their grounding is important and may be seen as the key to understand the foundations of Mathematics as well as Information Technology and especially functional programming. 

Despite a long research, the notion of functional is still a challenge. What does it mean to operate on higher order notions in a meaningful and effective way?  Actually, since functionals are abstract notions, it is the question about the computational grounding of the functionals.  

The first research related to computable functionals is due to 
Stefan Banach and Stanislaw Mazur  (1937) \cite{BanachMazur1937} \cite{Mazur1963}. 
The first definition of the functional related to primitive recursion schema of second order was introduced by Rosa P\'{e}ter (1951) \cite{Peter1} \cite{Peter2}. 

A significant contribution was done  by Andrzej Grzegorczyk (1955)  \cite{Grzegorczyk1955a} \cite{Grzegorczyk1955b} \cite{Grzegorczyk1957}, and by Daniel Lacombe \cite{Lacombe1955}.  
The paper by Grzegorczyk (1964) \cite{Grzegorczyk1964} is of particular importance. It was the first attempt to grasp the meaning of recursive (constructible) objects in all finite types.  

The first general approach to effective computability on objects on higher types was done by Stephen Cole Kleene (1959) \cite{Kleene1959b} and   \cite{Kleene1963}. Functionals were interpreted (coded) as natural numbers on which operations can be made. In this way the internal structure (construction) of functionals was lost. Ambiguity in treating natural numbers as functions (functionals) and natural numbers at the same time resulted in partial functionals, that is, not defined for some arguments. From computational point of view, partiality is the same as meaningless computation.    

The next approaches done by  Kleene (1959) \cite{Kleene1959a} (countable functionals)  and Georg Kreisel (1959) \cite{Kreisel1959} (continuous functional) were based on the intuition that any finite initial part of the output (value) of a functionals must be determined by a finite part of the corresponding argument (input). The functionals were supposed to be general (total), i.e. defined on the whole domain. The limit of these increased finite parts defines the functional in question. The very similar idea was employed by Dana Scott \cite{Scott1970} to define his famous Domain, and independently also by Yuri Ershov \cite{Ershov1972} \cite{Ershov1973b}. 

Partial functionals are more easy to formalize than total functionals in the very similar way as partial recursive functions versus general recursive functions. 

In Scott-Ershov Domain, effective functionals are defined as elements of the Domain that are of the form 
$ x = \bigsqcup_{n=1}^\infty e_n $, where  $\{e_1, e_2, ... , e_n, ... \}$ is a recursively enumerable subset of effective base (say $E$) of the Domain. The set $E$ is also a recursively enumerable subset of the Domain.  Since this definition involves the partial recursive functions, also the resulting functionals are partial, i.e. they may be not total. A partial recursive functions determining the set $\{e_1, e_2,  ... , e_n, ... \}$ does not reflect the structure of the functional it defines as well as the operation of the functional. 
Hence, it is rather an abstract definition. 

Based on Richard Platek PhD thesis (1966)  \cite{Platek1966}), Scott (1969) proposed (in his famous manuscript published in 1993 as \cite{Scott1993}) Logic of Computable Functions (LCF) that is a formal theory (interpreted in the Domain) for description and reasoning on computable functionals. The theory is also called PCF (for partial computable functionals), and also known under the name  Programming Computable Functions, see Robin Milner (1977) \cite{milner1977fully}, see also  Longley \cite{Longley05notionsof}.  Besides the fix-point operator and lambda operator, it contains operator  
 ``if'' roughly corresponding to if-then-else constructor in programming. 
 It turned out that this theory is not sufficient to define (as terms) all effective functionals defined in the Domain. Two additional operators are needed, that is, 
``por'' (parallel or), and ``exists'' introduced independently by  
V.Yu. Sazonov 1976 \cite{Sazonov1976} and Gordon D. Plotkin 1977 \cite{Plotkin1977}. The second operator corresponds to minimization $\mu$-operator from the definition of partial recursive functions. The extension is called PCF++. 

Note that the notion of functionals, as terms in PCF++ and their abstract mathematical interpretation as effective functionals in the Domain, concerns the partial functionals. 
However, partiality is the same as meaningless computation. From the (practical) computational point of view only total (general) functionals are essential. 
It seems that the proper understanding of the grounding of the notion of general functionals is still a challenge. 

For the comprehensive overview of the computable functionals, see  John Longley and  Dag Normann 2015 \cite{longley2015higher}, 
and  \cite{Longley05notionsof} on www.  
For a historical context of the research on functionals see  Robert Soare 1999 \cite{Soare99thehistory}.

The presented work is a continuation of the idea of Grzegorczyk \cite{Grzegorczyk1964} (who was inspired by the System T of Kurt  G\"{o}del \cite{Godel-Dialectica}) concerning recursive objects of all finite types. 

The approach presented in the current paper aims at the grounding of the functionals in hardware (actually in finite parametrized structures). This corresponds to the intuitionism of Luitzen E. J. Brouwer \cite{Brouwer1913} about the Foundation of Mathematics. 

%%%%%%%
%%%%%%%%%%%%%%
\section{Functionals and type theories} 

Actually, PCF++ is a type theory. However, there are terms in PCF++ (interpreted in the Domain) that correspond to non-terminating functions. It means that there are terms that have no normal form, if the term rewriting (as computation) is considered. In the sequel we are going to review shortly subsequent type theories (with dependent types) that have the normalization property. Higher order terms of these theories correspond to total functionals.  

The phrase {\em effectively constructed objects} may be seen as a generalization of the notion of {\em recursive objects} and general computable functionals. %Objects can be represented as finite (usually parameterized) structures. 
Universe of effectively constructible objects is understood here as an abstract collection of all generic constructible objects. 

In the Universe, constructability is understood literally, i.e. it is not definability, like general recursive functions (according to G\"{o}del-Herbrand) that are defined by equations in Peano Arithmetic along with  proofs that the functions are total, that is, defined for all their arguments. Objects are not regarded as terms in lambda calculus, and in combinatory logic. 

Most theories formalizing the notion of effective constructability (earlier it was computability) are based on the lambda abstraction introduces by Alonzo Church that in principle was to capture the notion of function and computation. Having a term with a free variable, in order to denote a function (by this term) the  lambda operator is applied. Unlimited application of lambda abstraction results in contradiction (is meaningless), i.e. some terms cannot be reduced to normal form. This very reduction is regarded as computation. Introduction of types and restricting lambda abstraction only to typed variables results in a very simple type theory.  

Inspired by System T, Jean-Yves Girard created system F \cite{Girard71}, \cite{Girard1989}; independently also by John C. Reynolds \cite{Reynolds}. Since System F uses lambda and Lambda abstraction (variables run over types as objects), the terms are not explicit constructions. System F is very smart in its form, however, it is still a formal theory with term reduction as computation; it has strong normalization property.  

Per Martin-L\"{o}f Type Theory (ML TT for short) \cite{Martin-Lof} was intended to be an alternative foundation of Mathematics based on constructivism asserting  that to construct a mathematical object is the same as to prove that it exists. This is very close to the Curry-Howard correspondence {\em propositions as  types}.  In ML TT, there are  types for equality, and a cumulative hierarchy of universes.  However, ML TT is a formal theory, and it uses lambda abstraction. Searching for a grounding (concrete semantics) for ML TT by the Author long time ago, was the primary inspiration for the Universe presented in this work.  

{\em Calculus of Inductive Constructions (CoIC)},  created by Thierry Coquand  and G\'erard Huet \cite{Coquand1986} and \cite{Coquand1988}, is a lambda calculus with a rich type system like the System F. It was designed for {\em Coq Proof Assistant} \cite{Coq}, and  can serve as both a  functional programming language and as a  constructive foundation for Mathematics. Agda is a dependently typed functional programming language based also on ML TT; it is also a proof assistant, see at www \url{wiki.portal.chalmers.se/agda/} . Like other functional programing languages, computation on higher order objects (functionals) is reduced to lazy evaluation, that is, to reducing a term denoting a functional to its normal form only if it is necessary.  

ML TT,  System F, and  CoIC are based on lambda and Lambda abstraction, so that in their syntactic form they correspond to the term rewriting systems. Functionals are defined there as terms. The problem is with their interpretation (grounding). The common view is that such grounding is not needed if the term rewriting is considered as computation. 

%The current paradigm in Computer Science states that the computation on higher order objects can be done only in the symbolic way, i.e. higher order objects (functionals) can be represented only as terms whereas computations on them can be done only  by term rewriting.  

%In this work the symbolic computation is challenged. It is an attempt to show that the syntactic constructions (i.e. terms  in System F, ML TT etc.), can have  explicit and concrete grounding as constructions. In this sense, it follows the idea of Grzegorczyk's combinators \cite{Grzegorczyk1964}, and in some sense also combinators of Haskell B. Curry \cite{Curry-cl} combinatory logic. %Note that combinatory logic is the theoretical basis for functional programming language Haskell. 

%Effective  construction of an object cannot use actual infinity. If it is an inductive construction, then the induction parameter must be shown explicitly in the construction. For any fixed value of the parameter the construction must always be a finite structure. The Universe presented in this paper is supposed to consist only of such objects. 
%Objects are not identified with terms whereas computations are not term rewritings.  

However, this view was challenged recently by  Homotopy Type Theory: Univalent Foundations of Mathematics (HoTT)(2013) \cite{HoTT}; a formal theory based on ML TT and CoIC. HoTT aspires to be another foundation for Mathematics alternative to set theory (ZFC), by encoding general mathematical notions in terms of homotopy types. According to Vladimir Voevodsky \cite{vv}(one of the creators of HoTT) the univalent foundations  are adequate for  human reasoning as well as for computer verification of this reasoning.
 Generally, any such foundations should consist of three components. 
The first component is a formal deduction system (a language and rules); for HoTT it is CoIC.  
The second component is a structure that provides a meaning to the sentences of this language in  terms of mental objects intuitively comprehensible to humans; for HoTT it is interpretation of sentences of CoIC as univalent models related to homotopy types. 
The third component is  a way  that enables humans to encode mathematical ideas in terms of the objects directly associated with the language.

The above phrases: {\em mental objects} and {\em mathematical ideas} are not clear. Actually, in the univalent foundations, these mental objects (as homotopy types) are yet another formal theory. It seems that the main problem here is the lack of a  grounding (concrete semantics) of these mental objects and mathematical ideas. The concept of equality (relation) plays extremely important role in ML TT and HoTT. However, a formal axiomatic description of the notion of equality of two object of the same type, and then higher order equality is not sufficient to comprehend the essence of the notion of equality and in general of the notion of relation. 

The origins of HoTT related to the notion of Continuum are interesting and are discussed in the companion paper \cite{C}.    

%The proposed Universe is not yet another formal theory of types. It is intended to be a grounding for some formal theories as well as a generic method for constructing objects corresponding to data structures in programming. 

%It seems that the same idea was investigated at least since the beginning of the XX century.  However, it was done in formal ways by Church lambda calculus, Curry combinatory logic, G\"{o}del System T, Grzegorczyk System, Martin Lof TT,  Girard System F, and Coquand CoIC to mention only the most prominent works.  The Universe is an attempt to understand these formal theories, that is, to comprehend their grounding (semantics). 

%Universe and the notion of higher order objects is strongly related to notion of  computable functionals (Stephen C. Kleene \cite{Kleene1959a}\cite{Kleene1959b}\cite{Kleene1963}, Georg Kreisler \cite{Kreisel1959}, Grzegorczyk  \cite{Grzegorczyk1955a} and   \cite{Grzegorczyk1955b}, as well as to Richard Platek $\&$ Dana Scott $PCF^{++}$ \cite{Scott69,Scott1993}) and to Scott Domain \cite{Scott1970} as a mathematical semantics of functional programing languages. 

% potrzebny jakis łącznik pomiedzy tym co wyżej a tym co niżej 

Functionals are also important for functional programming as well as for  computing hardware design. 

\section{Functionals as hardware }
\label{hof}

The current paradigm in functional programming states that the computations on higher order objects (functions and functionals in general) consist in symbolic manipulations, i.e. term rewriting according to syntactic rules.

Although symbolic computations make sense (like algebraic calculations), and the computations on higher order objects can be done (e.g. via recursion  equations), the intuition behind the functionals is that they are {\em ``objects''} that can be construed as concrete structures. 
 
The John Backus's (1977) \cite{Backus} original idea of function-level programming language was based on {\em ``programs as mathematical objects''}, where the objects (functions and functionals) are used directly in computations.
 
What are the higher order objects? Are they terms? Perhaps terms merely denote the objects?

The functional programing languages based on term rewiring (like Haskell and F\#) are considered as non-von Neumann programming languages. However, where is a corresponding non-von Neumann computer architecture? Perhaps John Backus was not right. Perhaps the notion of functional %(created by human intellect) 
is still far from being understood. %However, it is remarkable that human brain is not built according to the von Neumann architecture. 
 
If the terms along with rewriting rules are functionals, then functional programming may be seen still as {\em ``value-level programming''}, i.e. terms are of data type String. So that, term rewriting is an additional string processing (actually, it is byte processing) done according to the von Neumann computer architecture.  Although the symbolic computation is at higher level of abstraction, it still does not impose the essential change on the computer architecture. 

%Can functionals be realized as hardware? 
The idea of functionals as hardware is not new, see Mary Sheeran (1984) \cite{Sheeran1984}. 
The approach proposed by C$\lambda$aSH \url{http://www.clash-lang.org/} to realize higher order functionals is interesting. It goes from Haskell and its high-level descriptions (syntax) and via term rewriting (lazy evaluation as semantics) to a standard HDL (Hardware Description Language). However, this approach is still not sufficient for the paradigm shift. For a survey of functional HDL, see  Peter Gammie (2013)  \cite{gammie}. Hence, the existing approaches to {\em ``functionals as hardware''} explore the notion of functional based on symbolic computation paradigm, i.e. term rewriting as computation. However, in these approaches unlimited recursion may lead to non-terminating computations. So that, usually, these approaches are limited to higher order primitive recursion only, see Coq \cite{Coq} for example. 

It seems that the hardware technology is very close to shift the paradigm. A Coarse Grain Reconfigurable Architecture (CGRA) Bjorn De Sutter et al. (2013) \cite{de2013coarse}, is a type of processor architecture that can be reconfigured at runtime. The reconfiguration is done on the Function Unit (FU) level. That is, for an application, array of interconnected FUs can be dynamically configured. 
Here, functionals may be envisioned as generic mechanisms for the management of connections in the reconfigurable arrays. 
Hence, the concept of higher order computations as dynamic configuration of connections between hardware functional units is worth to be explored. 

Reconfigurable computing architectures (Russell Tessier et al. (2015)  \cite{tessier2015reconfigurable}), and reconfigurable system  (James Lyke et al. (2015)  \cite{lyke2015introduction}) are active research subjects. However, the correspondence between functionals and generic mechanisms for dynamic reconfigurations is still not recognized in the hardware design. 

 % Perhaps the most advanced approach was done by Dan Ghica et al. 2011 \cite{ghica2011geometry} that is limited, however, to compiling affine recursive programs.    

%Software and hardware are two different worlds. The first world consists of tightly coupled software (von Neumann programming languages) and extremely simple and ingenious von Neumann computing architecture. Computation is done on bytes sequentially on a single processor. 
%The second world is hardware with its rich diversity of possible computing architectures, and possible grounding for higher order computations on functionals. 

%A new concept of programming (different form the von Neumann languages) is needed to break the paradigm. Higher order computation can be envisioned as dynamic creation and reconfiguration of links between elementary circuits representing first order functions. 

% Perhaps the answer is in technology for developing large dynamically reconfigurable arrays of integrated circuits representing first order functions.  

There are also new trends in CGRA: data-flow graphs (Anja Niedermeier et al. (2014) \cite{niedermeier2014dataflow}, and Francesca Palumbo et al. (2016) \cite{palumbo2016dataflow}), full pipelining and dynamic composition  (Jason Cong et al. (2014) \cite{cong2014fully}), and overlay architecture  (Davor Capalija et al. (2014) \cite{capalija2014tile}, Abhishek Jain et al. (2016)  \cite{jain2016adapting}, and David Andrews et al. (2016) \cite{ma2016just}).

According to these new trends, functional units (FUs) as elementary first order functions are built on fine grained integrated circuits (e.g. ASICs). FUs form an overlay, if they are collected in an array, and connections between them are reconfigurable. Algorithm (to be realized in hardware) is modeled as a data-flow graph where the nodes correspond to FUs. If the array is sufficient rich in FUs and possible connections, the graph can be mapped into the array. If the graphs are acyclic, then the data flow if fully pipelined. Full pipelining is the best choice for hardware realization. Then, only appropriate buffering of input and output data of FUs is needed. For cycles, data-flow control and synchronization are necessary.  
 
Fully pipelined data-flows (as directed acyclic graphs) correspond to simple algorithms. Sophisticated algorithms use recursion and {\em ``while''} loops that enforce cycles in the graphs. The cycles can be eliminated if the graphs are dynamic in the sense that during execution, for a concrete value of the recursion parameter, the recursion node (or subgraph) can be unfolded to a acyclic subgraph to form new acyclic graph. This presupposes dynamic unfolding and edge configuration, i.e.  dynamic graph transformation, during execution. Recently, see  Mukherjee et al. (2017) \cite{Mukherjee}, this idea was explored   in a tool for automatic generation of directed acyclic dynamic graphs for sequential C programs. 

The transformation corresponds to the generic mechanisms for configuration of connections in large arrays of FUs.  
If the number of the connections is at most dozens (for simple computations), then dedicated mechanisms may be designed. However, if hundreds, thousands, and even more connections and FUs are needed, then the mechanisms must be generic, and must be based on higher order abstractions, i.e. higher order types and higher order functionals. 

What are these generic mechanisms (functionals)? How can they be grounded in hardware? From abstract mathematical point of view, the mechanisms correspond to transformations of acyclic directed graphs. It is a bit surprising, that if the graphs are considered as {\em ``terms''}, then the transformations correspond again to the {\em``term rewriting''}. However, the crucial point here is the hardware interpretation of the graphs and of the transformations. 

A preliminary framework exploring this idea is presented in the consecutive sections. 
The notion of function as well as higher order objects (functionals) is based on the elementary notions of type, object of type, type constructors, type of function, application of functional to object of higher type, and composition of two functionals. We show that these very elementary notions can be realized as hardware. Then, primitive type constructors, the primitive type of natural numbers, and primitive functionals are introduced. Higher order recursion schemata are constructed as functionals, i.e. as directed acyclic dynamic graphs. 
%
%In the full version of the paper, 
% to na konferencję 
Also dependent types, and relations on natural numbers are introduced. Finally, an example of pure function-level programming is presented. 
Hardware interpretation of the all introduced higher order objects and types is presented.  

The proposed framework is quite powerful, and its hardware grounding is important.  From the computational point of view, it seems that the proposed framework is reacher than Coq \cite{Coq} that is based on the symbolic computation paradigm.

%%%%
%%%%%%%%%
\section{Foundations }
\label{gniazda-wtyczki}
 
Let $A$, $C$ and $B$ denote primitive data types (e.g. Integers). 
Let us consider two first order functions $f: A \rightarrow B$, and $g: B \rightarrow C$.  If the functions are interpreted as integrated circuits, then  application of function $f$ to $a$ of type $A$, usually denoted by $f(a)$, is clear, and means that the signal corresponding to $a$ is at the input of $f$. Then, the signal is processed by $f$, and the result (also as a signal) is at the output of $f$. The composition of the two function, i.e. $f\circ g$  such that   $(f\circ g) (a) = g(f(a)) $, consists in connecting the output of $f$ to the input of $g$. Note that $f\circ g$ is a function of type  $A \rightarrow C$. 

Generally, first order function consists of sockets (corresponding to input data), body (where the data is processed), and plugs (corresponding to output) for temporal storing the results of processing. 

\begin{figure}[h]
	\centering
	\includegraphics[width=0.8\textwidth]{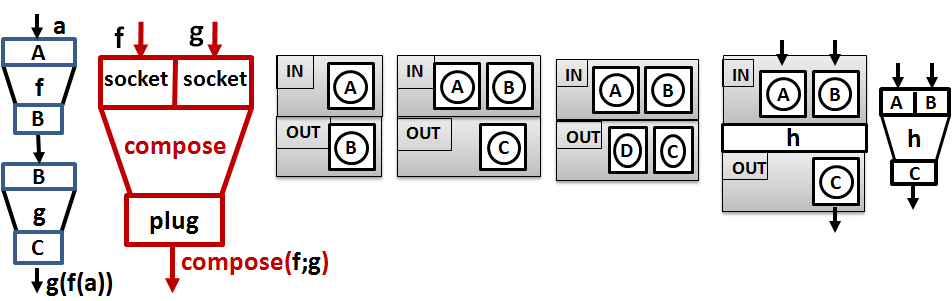}
	\caption{On the left, amorphous function composition, and function composition as a functional. In the middle, types as boards of sockets and plugs: $(A \rightarrow B)$, and $((A; B) \rightarrow C)$, and $((A; B) \rightarrow (D;C))$. On the right, two pictorial representations of function $h$ of type $((A; B) \rightarrow C)$ }
	\label{rys1-2}
\end{figure}

Function $f$ has input (socket) of type $A$ and output (plug) of type $B$. Function $g$ has input (socket) of type $B$ and output (plug) of type $C$.  Since the plug of $f$ is of the same type as the type of socket of $g$, the directed connection between the plug and the socket (putting plug into socket) means the composition of functions $f$ and $g$, see Fig.\ref{rys1-2}. %Function composition is one of the basic notions in Mathematics. 

For the first order functions, application and composition have the clear hardware interpretation. We are going to show that also higher order application and composition can be grounded in hardware. In order to do so, the higher order types must be interpreted in hardware. 
Type of a first order function is a board consisting of two parts. The first part is for sockets, whereas the second part is for plugs, see Fig. \ref{rys1-2}. Since the sockets and the plugs are of primitive types, the hardware interpretation of the board is simple. 

\begin{figure}[h]
	\centering
	\includegraphics[width=0.8\textwidth]{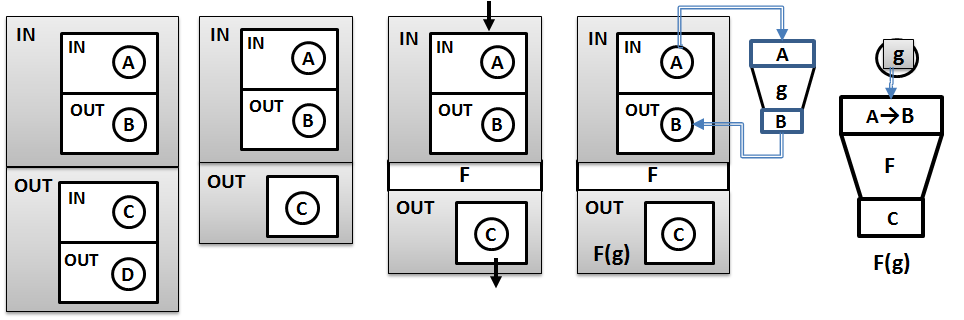}
	\caption{On the left, higher order types: $(A \rightarrow B) \rightarrow (C \rightarrow D)$, and  $(A \rightarrow B) \rightarrow C$. In the middle, amorphous higher order application of functional $F$ to function  $g: A\rightarrow B$. The result $F(g)$ is at the plug of the functional. On the right, a pictorial representation of the application }
	\label{rys3-4}
\end{figure}

Higher order types are construed as nested boards of sockets and plugs of primitive types. The nested board for the type $(A \rightarrow B) \rightarrow C$ (see Fig. \ref{rys3-4}), consists of a socket of type  $A \rightarrow B$, and a plug of type $C$. The socket $A \rightarrow B$ itself is a sub-board consisting of a sub-socket of type $A$, and a sub-plug of type $B$.   
The nested board for type  
$(A \rightarrow B) \rightarrow (C \rightarrow D)$, see Fig. \ref{rys3-4}, consists of one socket of type $A \rightarrow B$, and one plug of type $C \rightarrow D$.  
The boards for more complex types can be constructed in an analogous way.

{\bf Higher order amorphous application.}  
For a functional $F$ of type $(A \rightarrow B) \rightarrow C$, its application to function $g: A \rightarrow B$ is realized in the simple way, see Fig.  \ref{rys3-4}. 
The socket of the functional $F$ is of type $A\rightarrow B$. The application consists in establishing appropriate connections (directed links) between the socket of $F$ and the socket and plug of the function $g$. The directed links correspond to the data flow. The first link is between sub-socket of the socket of functional $F$ and the socket of $g$. The second link is between the plug of $g$ and the sub-plug of the socket of $F$.   The result $F(g)$ is of type $C$, and is at the plug of the functional.

\begin{figure}[h]
	\centering
	\includegraphics[width=0.7\textwidth]{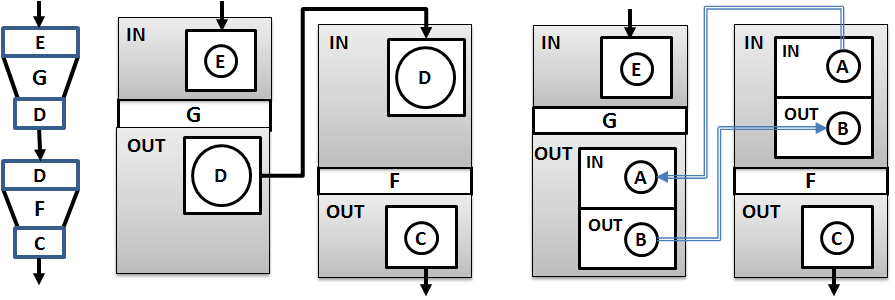}
	\caption{$D$ denotes the type $A \rightarrow B$. Amorphous composition of  functional $F$ of type $E \rightarrow D $ and functional  $G$ of type $D\rightarrow C$ }
	\label{acomposition}
\end{figure}

{\bf Higher order amorphous composition.} 
Let $D$ denote the type $A \rightarrow B$. For a functional $G$ of type $E \rightarrow D$ and a functional $F$ of type $D \rightarrow C$, their amorphous composition $G\circ F$ is shown in Fig. \ref{acomposition}. It consists in making a connection between the plug of functional $G$, and the socket of functional $F$. Since the type of the plug and of the socket is of  higher order type, (i.e.  $A \rightarrow B$),  this {\em higher order connection} consists of two first order connections, see Fig.  \ref{acomposition}. The first connection is between the sub-socket of the socket of $F$ and the sub-socket of the plug of $G$, whereas the second connection is between the sub-plug of the plug of $G$ and the sub-plug of the socket of $F$.  
Note that for $e$ of type $E$, \ \ $G(e)$ is a function of type $A  \rightarrow B$, and is at the plug of $G$. Denote this function by $g$.  Higher order amorphous application (see Fig. \ref{rys3-4}) of $F$ to $g$, that is, $F(g)$ corresponds exactly to the connections in Fig. \ref{acomposition}. It means that  $F(G(e))$ is the same as $F(g)$, and the same as $(G\circ F)(e)$. 

Hence, the above hardware interpretations of the higher order amorphous application and higher order amorphous composition make sense. Note that if the types $A$, $B$, and $C$ are of higher order, then in the above hardware interpretations, the first order connections become higher order connections. So that, the interpretations are generic. 

Functional versions of higher order composition, and higher order application are presented below. 
 \begin{figure}[h]
	\centering
	\includegraphics[width=0.87\textwidth]{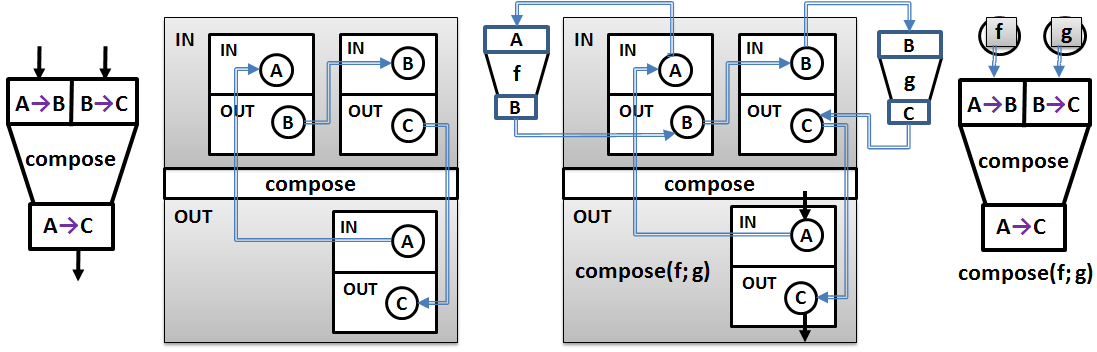}
	\caption{ Functional $compose_{A,B,C}$, and application of the functional to the two functions: $f$ and $g$. The result is at the plug of $compose$. On the right, a pictorial representation of the application }
	\label{nc-5}
\end{figure} 
Composition as a functional, i.e. $compose_{A,B,C}$ of type $((A \rightarrow B); (B \rightarrow C)) \rightarrow (A \rightarrow C)$ for composition of two functions, $f$ of type $A\rightarrow B$, and $g$ of type $B\rightarrow C$, can be constructed by making appropriate connections in a nested board of sockets and plugs, see Fig. \ref{nc-5}. Application of $compose_{A,B,C}$ to functions $f$ and $g$ results in their composition, i.e. $compose_{A,B,C}(f;g)$ is the same as $f\circ g$. To grasp the idea just follow the directed links in the Fig. \ref{nc-5} (second picture from the right) from sub-socket of type $A$ of the plug of type $A\rightarrow C$ of the functional via sockets and plugs to the sub-plug of type $C$ of the plug of the functional.  

\begin{figure}[h]
	\centering
	\includegraphics[width=0.73\textwidth]{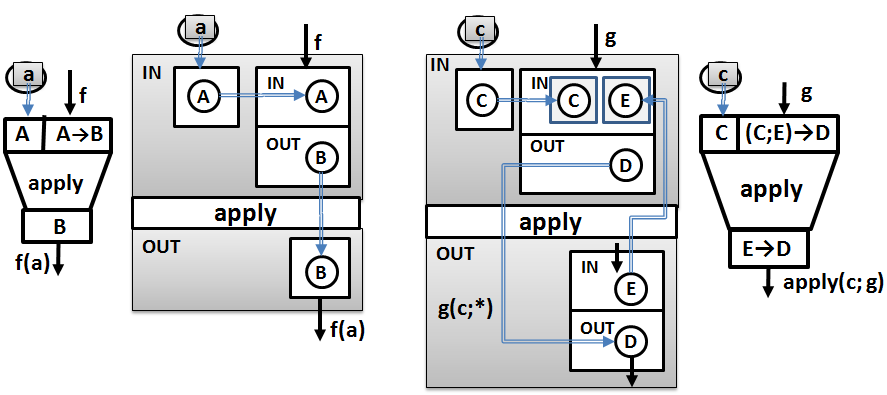}
	\caption{Simple application $apply_{A\rightarrow B, A}$ of type  $(A; (A\rightarrow B)) \rightarrow  B$, and more complex application of type $(C;((C;E) \rightarrow D))\rightarrow (E \rightarrow  D)$ }
	\label{rys7}
\end{figure}

Application as a functional, i.e. $apply_{A\rightarrow B, A}$ is of type  $(A; (A\rightarrow B)) \rightarrow  B$ such that for any $f$ of type $A\rightarrow B$, and any $a$ of type $A$, \ \ $apply_{A\rightarrow B, A}(a; f)$ is $f(a)$. The construction of the functional is shown in Fig. \ref{rys7}. 

Note that if the types $A$, $B$, $C$ and $D$ are of higher order, then the corresponding connections are of higher order. Hence, the above hardware interpretations are generic.  

The conclusion is that higher order types may be realized as hardware in the form of nested boards of sockets and plugs of the primitive types. Higher order applications and compositions may be realized as hardware by making appropriate connections in the nested boards. 

The above (a bit informal) framework it simple.  However, as the foundation,  it is powerful enough to construct higher order primitive recursion schemata, and in general (it seems) the intuitionistic second order Arithmetics. 
The framework is the necessary basis for the next sections. 

%%%%%% 
%%%%%%%%%%%%
\section {Level zero } 

We are going to construct a Universe (of construable objects) as the hierarchy of universes starting with the level zero. This is similar to the universes in 
Per Martin-L\"{o}f Type Theory (ML TT for short) \cite{Martin-Lof}, and 
 Calculus of Inductive Constructions (CoIC)  created by Thierry Coquand  and G\'erard Huet \cite{Coquand1986} \cite{Coquand1988}. However, ML TT and CoIC are merely formal type theories with semantics based on term rewriting. 
At the level zero of the Universe, primitive type constructors, a primitive type, and corresponding primitive operations are introduced.  

%%%%%%  
%\subsection{ Notations } 

Let $a:A$ denote that object $a$ is of type $A$. The word {\em ``operation''} will be used instead of {\em ``function''} and {\em ``functional''}. 
Type of operation is denoted by $A^s \rightarrow B^p$ where $A^s$ denotes a socket of type $A$, whereas $B^p$ denotes a plug of type $B$.  
The general form of operation with multiple sockets and multiple plugs is \\
$
g: \ (A^s_1; A^s_2; ... ; A^s_k) \rightarrow  (B^p_1; B^p_2; ... ; B^p_l)
$. 
Upper indexes denoting sockets and plugs will be frequently omitted. The order of sockets (plugs) is not relevant. 

Amorphous application of operation  $f:A \rightarrow B$ to object $a:A$ is denoted by $f(a)$. For operations with multiple input, application may be partial, e.g. application of $g$ only to $a_i:A_i$ and $a_j:A_j$ is denoted by $g(a_j; a_i; *)$.

%%%%%%%%
\subsection{ Simple type constructors }
\label{type-constructors}

\begin{figure}[h]
	\centering
	\includegraphics[width=0.9\textwidth]{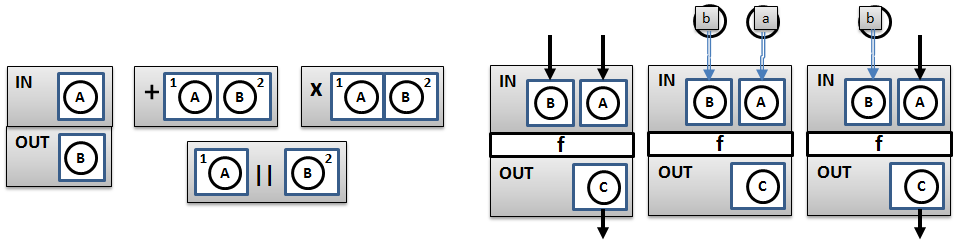}
	\caption{ On the left, a function type, disjoin union, product, and exclusive disjunction. On the right, applications of $f:(B^s;A^s) \rightarrow  C^p$ to two objects, and to one object }
	\label{rys5-6}
\end{figure}

Having in mind the hardware interpretation of types as boards of sockets and plugs, there are four basic type constructors. 
Let $A$ and $B$ denote types. 
\begin{itemize}
\item  $\times$, product of two types denoted by $A\times B$; an object of this type consists of two objects, one of type $A$ and one of type $B$. 
\item  $+$, disjoin union $A+B$; object of this type consists of either of an object of type $A$ and the pointer to type $A$, or it consists of an object of type $B$ and the pointer to type $B$.  
\item  $\rightarrow$, arrow as the constructor of function type $A \rightarrow B$; where $A$ is socket and $B$ is plug.  In its general form there may be a finite number of sockets, and a finite number of plugs. Some of the types of the sockets as well as the plugs may be the same.   
\item $||$,  exclusive disjunction $A||B$; special type used for plugs. It means either $A$ or $B$. It is used for the plug of primitive destructor operation $get_{A,B}$ introduced in the next subsection, and for the plug of primitive operation {\tt if-then-else} introduced in section \ref{if-then-else}. 
\end{itemize}
Given two types (boards) $A$ and $B$, the constructors produce a third board as shown in Fig. \ref{rys5-6}.  The boards  $A$ and $B$ have their identifiers in the resulting board.    

For operation $g:A^s \rightarrow C^p$, putting object $a:A$ into socket  $A^s$ means amorphous application with the result $g(a)$ available at plug $C^p$. 
For operation $f:(B^s; A^s) \rightarrow  C^p$, the application result $f(b,a)$ is at the plug $C^p$, see Fig. \ref{rys5-6}. However, a partial application $f(b, *)$  is an operation of type  $A^s \rightarrow  C^p$.

%%%%%%%
\subsection{ Object constructors and destructors } 
\label{k-d}

For $a:A$ and $b:B$, object constructors are introduced in the following way:  
\begin{itemize}
\item for product, $join_{A,B}$ is a primitive operation of type $(A;B)\rightarrow (A\times B)$ such that $join_{A,B}(a;b)$ is object of type $A\times B$ denoted as pair $(a,b)$.
\item for disjoin union, $plus^A_{A,B}: A \rightarrow (A+B)$  and  $plus^B_{A,B}: B\rightarrow (A+B)$ are primitive operations such that 
 $plus^A_{A,B}(a)$ and  $plus^B_{A,B}(b)$ are objects of type $A+B$, denoted respectively by $(1.a)$ and $(2.b)$, where $1$ is the pointer   to type $A$, and $2$ to $B$.  
\item for arrow there are two primitive object constructors:
\begin{itemize}
\item constant operations, $const_{A,B}: A \rightarrow (B \rightarrow A)$ such that for any $a:A$ operation  $const_{A,B}(a): B \rightarrow A$ is constant operation, i.e. for any $b:B$, \ \ $const_{A,B}(a)(b)$ is always $a$.  
\item identity operations,  $id_A:A\rightarrow A$ such that for any $a:A$, the result $id_A(a)$ is $a$. 
\end{itemize}
\end{itemize}

Object destructors are as follows. 
\begin{itemize}
\item $proj_{A,B}:(A\times B) \rightarrow (A; B)$. Note that this operation has two plugs. For $(a,b)$ of type $A\times B$, \ \ $proj_{A,B}((a,b))$ consists of two objects: $proj^A_{A,B}((a,b))$ that is $a$, and $proj^B_{A,B}((a,b))$ that is $b$.   
\item $get_{A,B}:(A+B) \rightarrow (A||B)$.  For  $(1.a)$ and $(2.b)$ of type $A+B$,\ \  $get_{A,B}((1.a))$ is $a$, and $get_{A,B}((2.b))$ is $b$. 
\item $apply_{A\rightarrow B, A}$ of type $((A \rightarrow B);A)\rightarrow B $. For any $f:A \rightarrow B$ and $a:A$,  \ \ $apply_{A\rightarrow B, A}(f;a)$ is $f(a)$. Note that amorphous application $()$ is also a destructor for arrow.    
\end{itemize} 

The constructors and destructors have simple hardware interpretations. The interpretations of application and partial application have already been  shown in Fig. \ref{rys7}.

%%%%%
\subsection{ Composition as operation }

Simple composition of operation $f:A^s\rightarrow B^p$ and operation $g:B^s\rightarrow C^p$ consists in establishing the connection between plug $B^p$ and socket $B^s$. This very establishing is the amorphous composition; the only requirement is that the type of the  plug and the type of the socket are the same. The amorphous composition is denoted (as usual) by the symbol $\circ$, and applied to the two operations $f$ and $g$ is denoted by $f \circ g$. 

There is also operational version of composition, i.e. for fixed types of two operations, if the type of a plug of one operation is the same as the type of a socket of the second operation, then the composition (making connection between the plug and the socket) can be done by special operation.   

Composition as operation $compose_{A,B,C}: ((A\rightarrow B);(B\rightarrow C)) \rightarrow  (A\rightarrow C)$ has already been constructed, see Fig. \ref{nc-5}.  
The most simple composition, i.e. $compose_{A,A,A}: ((A\rightarrow A);(A\rightarrow A)) \rightarrow  (A\rightarrow A)$ will be used in the  construction of primitive operation $Iter$ in the subsection \ref{Iter}.

%%%%%%%%%%
\subsection{ Operation $Copy$ }

Given an object, the operation $Copy$ produces a copy of the object. Hardware interpretation of the operation is neither simple nor obvious. One of possible interpretation assumes that every object contains the description of its construction.  

For already constructed object $a$, $Copy(a)$ returns two objects. The fist one is denoted by $Copy^1(a)$; it is the original object $a$. The second object denoted by $Copy^2(a)$ is a copy of $a$. Frequently, the copy will be denoted by $a'$.
$Copy$ is an amorphous operation. However, in constructions, its typed version can be used as operation $Copy_A: A\rightarrow (A;A)$. 
Once an object is used in a construction, it can not be used again. In order to distinguish between the original and its copies, the following notation is used. Symbol without apostrophe denotes original object. Symbol with one apostrophe (several apostrophes) denotes copy (consecutive copy).  For example,  $A'$ and $A''$ are copies of $A$ whereas $a'$, $a''$ and $a'''$ are copies of $a$.  Sometimes apostrophes will be omitted.    
%So far no primitive type was introduced. 

%%%%%%%%%
\subsection{The primitive type of natural numbers } 
\label{liczby_naturalne}
\begin{figure}[h]
	\centering
	\includegraphics[width=0.75\textwidth]{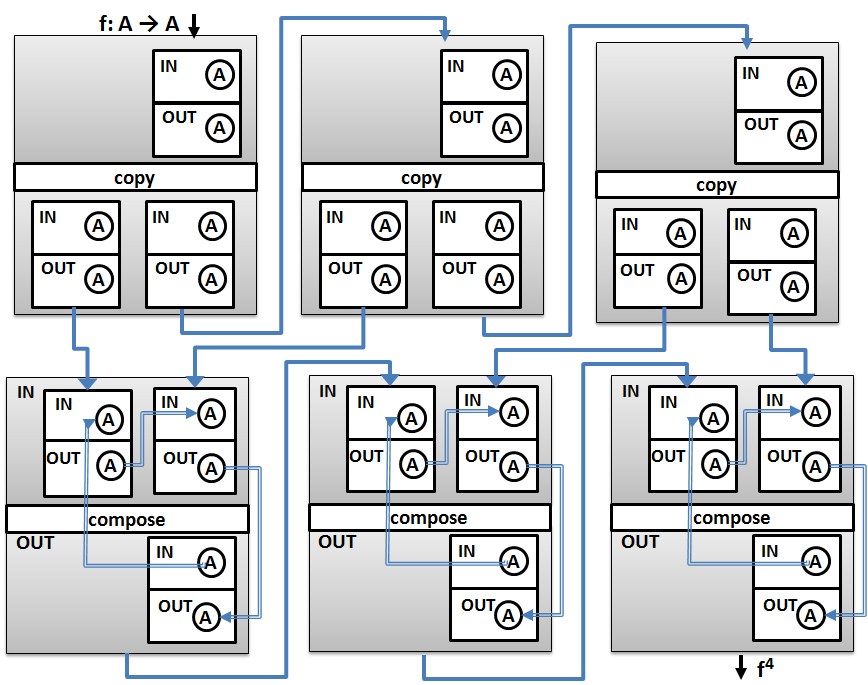}
	\caption{Operation $Iter(4;*)$ applied to $f$.  Note that every connection between plug of type $A\rightarrow A$ and  socket of type $A\rightarrow A$ is, in fact, a higher type connection }
	\label{rys5}
\end{figure} 
Copy the unit signal from the transmission channel, and the result put in the channel at the beginning.  Repeating this means natural numbers. Starting with a single unit signal (denoting number 1), the consecutive repetitions give next natural numbers, i.e. the first repetition results in two unit signals, the second one in three signals, and so on. This very repeating is the successor operation denoted by $Suc$. 
 
If the above intuition is applied to an operation of type $A\rightarrow A$ (here $A$ is an arbitrary type)) instead of the unit signal, then it is exactly the approach to define natural numbers proposed by Alonzo Church in  his lambda calculus, and also the one used in System F (Jean-Yves Girard (1971)  \cite{Girard71}, and John Reynolds (1974) \cite{Reynolds}). Natural number (say $n$) is identified with the amorphous iteration, i.e. it can be applied to any operation (with input and output of the same type), and returns $n$-times composition of this operation. 

Let us accept the first interpretation. Then the primitive operations, i.e.  successor $Suc$ and predecessor $Pred$, have natural interpretation. $Suc$ consists in coping the original unit signal and joining the result to what already has been done. $Pred$ is interpreted as removing from the channel the first unit signal, if the channel is not empty. 
Let $N$ denote the type of natural numbers.  

The type of natural numbers is the basis for constructing sophisticated objects. The constructions are inductive, that is, each construction depends  on a finite number of parameters. Although, potentially an inductive object is infinite, at any moment of time of the process of its construction, the values of the parameters are bounded, so that it is always a finite structure. 

%%%%%%%%%
\subsection{Iteration }
\label{Iter}
 
Iteration (indexed by type $A$) is the operation   
$Iter_A: (N;(A\rightarrow A)) \rightarrow  (A\rightarrow A)$, 
such that for any $n:N$ and any operation $f:A \rightarrow A$, the result  $Iter_A (n;f)$ is $n$-times composition of operation $f$. 
Parameter $n:N$ determines how many copies of $f$ and copies of operation  $compose_{A,A,A}$ must be used to produce the result. 

Construction of the result depends on $n$ and consists in dynamic configuration of connections between plugs and sockets as it is shown in Fig. \ref{rys5} for $Iter_A(4;*): (A\rightarrow A) \rightarrow  (A\rightarrow A)$.  

%%%%%%%%%%
\subsection{Operation $Change$ }

The next primitive operation is $Change_A$ of type  
$(N; A; (N \rightarrow A))\rightarrow (N \rightarrow A)$ such that for any $n:N$, object $a:A$, and operation $q:N \rightarrow A$, the result $Change_A(n;a; q)(n)$ is $a$. For $k$ different than $n$, result  $Change_A(n;a; q)(k)$ is $q(k)$. 
That is, for a sequence $q$ of objects of type $A$, and $a:A$, change $n$-th element of the sequence $q$ to $a$, i.e. $q(n)$ to $a$. 
Operation $Change$ corresponds to {\tt if-then} and {\tt case} statements in programming. 
For hardware interpretation of $Change_A$, evaluation of the condition {\em``If input is $n$''} must be realized. This corresponds to a primitive relation $Equal_N$ on type $N$ that will be introduced in section \ref{relacje}.

%%%%%%%%%%% UWAGA: tutaj jest istotna zmiana w stosunku do TO i dotyczy currying, jest to inna i naprawdę już operacja 
\subsection{ Currying and uncurrying}
\label{currying}

{\em Currying} is a syntactical rule to transform a term denoting a function with two or more variables to equivalent (nested) term with one outer variable; the other variables are hidden inside the term. It was introduced by Moses Sch\"{o}nfinkel in 1924 and later developed by Haskell Curry. 
\begin{figure}[h]
	\centering
	\includegraphics[width=0.5\textwidth]{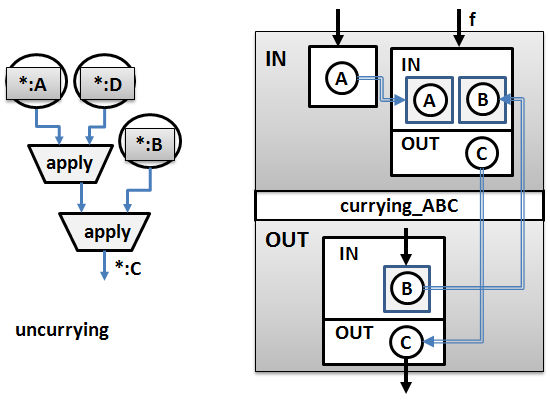}
	\caption{On the left, construction of operation $uncurrying_{A,B,C}$; here $D$ denotes type $A\rightarrow (B\rightarrow C)$. On the right, construction of operation $currying_{A,B,C}$}
	\label{rys8nowy}
\end{figure} 
Currying as well as {\em uncurrying} (i.e. the reverse transformation dual to currying) can be represented as operations.

Operation $f:(A;B) \rightarrow  C$ is transformed by $currying$ into  operation $g:A\rightarrow (B\rightarrow C)$ such that $f(a;b)$ is $g(a)(b)$. 

Operation $currying_{ABC}: (A; ((A;B)\rightarrow C)) \rightarrow (B\rightarrow C)$ is constructed in Fig. \ref{rys8nowy}. Actually, it is a board consisting of two sockets, one type  $A$ and the second one of type $(A;B)\rightarrow C)$, a plug of type  $B\rightarrow C$, and the appropriate connections.  Partial application of $currying_{ABC}$ to $f: (A;B)\rightarrow C$, i.e.  $currying_{ABC}(*;f)$ is operation $g:(A\rightarrow (B\rightarrow C)$ such that $f(a;b)$ is $g(a)(b)$.

Operation $uncurrying_{ABC}$ of type $(A; (A\rightarrow (B\rightarrow C));B) \rightarrow C)$ is constructed in the following way, see also Fig. \ref{rys8nowy}. 
Compose 
 $apply_{(A\rightarrow (B\rightarrow C)), A}: (A; (A\rightarrow (B\rightarrow C)) \rightarrow (B\rightarrow C) $ 
 and  $apply_{ (B\rightarrow C), B}: (B; (B\rightarrow C)) \rightarrow C$.  
It is done by connecting the plug of type $B\rightarrow C$ of the first operation to the socket of the type $B\rightarrow C$ of the second operation. The result is the required $uncurrying_{ABC}$. It has three sockets. By partial application (i.e. putting $g$ of type $A\rightarrow (B\rightarrow C)$ into the appropriate socket of $uncurrying_{ABC}$), we get operation $uncurrying_{ABC}(*;g;*)$ of type $(A ; B) \rightarrow C$. 

 Note that partial amorphous application is the crucial element of the above two constructions. 

%%%%%%
%%%%%%%%%%%%%
\section{The primitive recursion schema}
\label{rec}

The schema of primitive recursion for operations of the first order (from natural numbers into natural numbers) is clear. However, it is not so obvious for operations of higher types. The recursion schema for second order operations was introduced by R\'{o}zsa P\'{e}ter (1936) \cite{Peter}. 
G\"{o}del's System T (1958) \cite{Godel-Dialectica}, and  Andrzej Grzegorczyk`s System (1964)  \cite{Grzegorczyk1964} are based on the recursion on higher types. Grzegorczyk's iterators (as primitive recursion schemata indexed by types) are considered as objects.  
Haskell Curry (1964) \cite{Curry} defined Grzegorczyk's iterators as terms in combinatory logic using pure iteration combinator corresponding to the operation $Iter$ introduced in section \ref{Iter}. 
Jean-Yves Girard (1989) \cite{Girard1989} defined higher order primitive recursion schemata as terms in his System F. 

The higher-order recursion is still of interest mainly because of its application in programming. However, recent works are based on formal approaches. For the G\"{o}del-Herbrand style approach see Lawrence Paulson (1984) \cite{Paulson};  it is still not clear what is the meaning of equality for objects of higher types.
In Martin Hofmann (1999) \cite{hofmann1999semantical}, and Jiri Adamek et al. (2011)  \cite{Adamek}, a category-theoretic
semantics of higher-order recursion schemes was presented. In programming, see Ana Bove et al. (2005) \cite{Bove}, recursive algorithms are defined by equations in which the recursive calls do not guarantee a termination. 
Finally, Carsten Schurmann, 
Joelle Despeyroux, and  Frank Pfenning  (2001) \cite{Schurmann} propose an extension of the simply typed lambda-calculus with iteration and case constructs. Then, higher order primitive recursive functions  are expressible through a combination of these constructs. Actually, they defined terms that correspond to Grzegorczyk's iterators. A construction of the iterators as functionals (operations) are  presented below.   

%%%%%%%
\subsection{ Grzegorczyk's iterator  }
\label{G}

Although the Grzegorczyk System was intended to be constructive, it is still a formal theory. 
Grzegorczyk's iterator denoted here by  $R^A$ is a primitive (in Grzegorczyk System) term  of type 
$A \rightarrow ( (N\rightarrow (A\rightarrow A)) \rightarrow  (N\rightarrow A) )$ 
that satisfies the following equations: 

for any $a:A$, \ \ $c:N\rightarrow (A\rightarrow A)$, and  $k:N$

$R^A(a)(c)(1) = a $

$R^A(a)(c)(k+1) = c(k)(R^A(a)(c)(k))$
\\
The notation for nested applications is that  
$f(a)(c)(k)$ is the same as $((f(a))(c))(k)$. 

Note that in the above equations, the equality is for objects of type $A$. For the type of natural numbers it is clear, see section \ref{relacje}. However, it isn't so for higher order types. 
In a formal theory, the axioms of equality for higher order types must be introduced.  

By applying currying and uncurrying, $R^A$ can be interpreted equivalently as operation of type \\ 
$(N\rightarrow (A\rightarrow A)) \rightarrow  (N\rightarrow (A\rightarrow A))$ 
denoted by $\bar{R}^A$ such that $\bar{R}^A(c) (k) (a)$  is the same as  $R^A(a)(c)(k)$. Now, the definitional equations above can be rewritten as:

$\bar{R}^A(c)(1)(a) = a$ 

$\bar{R}^A(c)(k+1)(a) = c(k)(\bar{R}^A(c) (k) (a) )$ 
\\
In this new and equivalent form the iterator is the operation that for $c:N\rightarrow (A\rightarrow A)$, i.e. a sequence of operations of type $A \rightarrow A$,  
 produces object (sequence) 
 $\bar{c}: N \rightarrow (A\rightarrow A)$ such that $\bar{c}$ is $\bar{R}^A(c)$. 
 First element of this sequence, i.e.  $\bar{c}(1)$, is the identity operation on $A$ such that for all $a:A$, $id_A(a)$ is $a$.  
 
The element  $\bar{c}(k+1)$, that is $\bar{R}^A(c)(k+1)$, is the composition of operation $c(k)$ and operation $\bar{c}(k)$ that is $(\bar{R}^A)(c)(k)$. So that,  $\bar{c}(k+1)$ is the composition of the first $k$ elements of the sequence $c$.

In this equivalent form, the Grzegorczyk iterator is simple. However, its construction needs some effort.   

Girard's recursion operator (indexed by type $A$, however here denoted by $R$) is defined as a Lambda term in System F.   
Definition of $R$ is based on the interpretation of natural numbers as amorphous iterating operators. Applying number $n$ to arbitrary term (denoting a function with input and output of the same type) means to compose $n$-times the function with itself. 
The recursion operator $R$ is of type $A \rightarrow  ((A\rightarrow (N\rightarrow A)) \rightarrow  (N \rightarrow  A)) $, and has the following property. 
For any $a:A$, $v: A\rightarrow (N\rightarrow A)$ and $k:N$,

$R(a)(v)(1) = a$

$R(a)(v)(k + 1) = v( R(a)(v)(k) )(k)$
\\
The equalities above must be understood (according to Girard)  that both sides can be reduced by term rewriting to the same normal form. 

Apply currying and uncurrying to $R$ in the similar way as for the Grzegorczyk's iterator. 
Note that $(A\rightarrow (N\rightarrow A))$ and $A$ can be swapped in the type of $R$, so that we get operation of type
$(A\rightarrow (N\rightarrow A)) \rightarrow  (A \rightarrow  (N \rightarrow  A))$. 
Then, in the first and the second segment,  $N$ and  $A$ can be swapped, so that we get the operation of type   
$(N\rightarrow (A\rightarrow A)) \rightarrow  (N\rightarrow (A\rightarrow A))$
denoted by $\bar{R}$ such that 

$\bar{R}(\bar{v})(1)(a) = a$, where  $\bar{v}:N\rightarrow (A\rightarrow A)$ satisfies $\bar{v}(k)(a) = v(a)(k)$,  

$\bar{R}(\bar{v})(k + 1)(a) = \bar{v}(k)( R(\bar{v})(k)(a) )$
\\
In this form $\bar{R}$ is exactly the same as  Grzegorczyk's iterator, i.e. for a sequence of operations of type $ A\rightarrow A $ as input, it returns the output sequence where its $(n+1)$-th element is the composition of the first $n$ elements of the input.

However, this cannot be taken literally that the input, as a sequence, is taken as whole (as actual infinity) by the Grzegorczyk's iterator (and Girard's operator), and returns a complete sequence as its output. From the syntactical point of view it is acceptable as a definition of a term, however, not as a construction. Actually, the parameter $n:N$, that refers here to the $n$-th  element of the input sequence, must refer to the construction parameter. It will be clear in the following construction. 

%%%%%%%%%%% 
\subsection{Construction of Grzegorczyk's iterator } 

Let $C$ denote the type $(N\rightarrow (A\rightarrow A))$. 
We are going to construct operation $iterator: C \rightarrow  C$, such that for any input object (sequence)  $c:C$, the operation returns object (sequence)   $\bar{c}:C$, such that $n$-th element of $\bar{c}$ is the composition of the first $n$ elements of sequence  $c$. Although, it is almost the same as for Grzegorczyk's iterator, it will be clear that at any construction step,  $iterator$ is a finite structure, that is, the parameter $n:N$ in the construction of $iterator$ refers always and only to this finite structure, i.e. to the first $n$ elements of $c$ and of $\bar{c}$. 

The detailed and explicit construction of this operation is simple. It is worth to carefully analyze this construction in order to grasp the full meaning of the idea of the constructability proposed in the paper. 

\begin{figure}[h]
	\centering
	\includegraphics[width=0.77\textwidth]{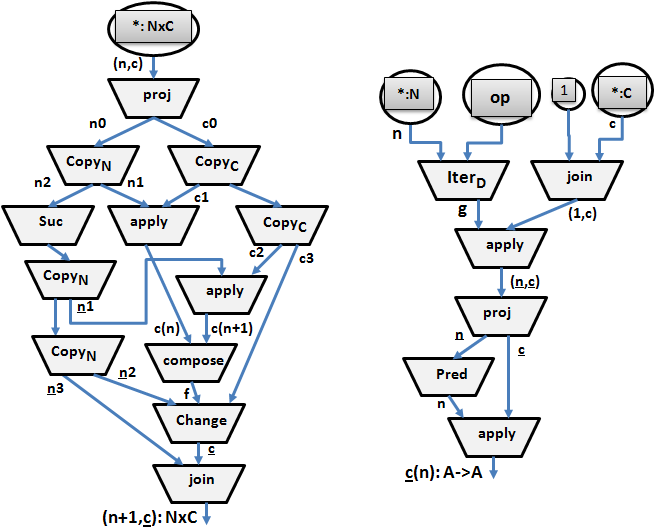}
	\caption{ Construction of operation $op: (N\times C) \rightarrow (N\times C)$, and operation $Rec_A:(N;C)\rightarrow (A \rightarrow A)$  }
	\label{rys6-7}
\end{figure}

The construction needs two auxiliary operations $op$ and $Rec_A$ constructed in  Fig. \ref{rys6-7} in the form of two acyclic directed graphs. The nodes of the graphs denote primitive operations or already constructed objects and operations. In the second graph, its two input (initial) nodes are labeled by $op$, and by $1$. The link between $op$ and operation $Iter_D$ corresponds to amorphous application of $Iter_D$ to $op$. Analogously for $1$ and operation $join$. The rest of the input nodes indicate the input types. The all links between internal nodes corresponds to amorphous composition, i.e. a connection between a plug of one operation to a socket of another operation; where the type of the plug and the type of the socket are the same. 
The graphs may be viewed as fully pipelined data-flows. The graph of operation $op$ is a simple data-flow.  The data-flow of the  graph of $Rec_A$ is nested and is unfolded if the parameter $n:N$ is instantiated to a concrete value. The unfolding is done by by applying $Iter_D$ to $n$, so that the $n$ copies of operation $op$ are composed.

The  operation $op$ takes pair $(n,c)$ as the input, and returns new pair $(n+1, \underline{c})$ where  $\underline{c}$ differs from $c$ only on the $(n+1)$-th element, that is, $\underline{c}(n+1)$ is the composition of two operations $c(n)$ and $c(n+1)$.  

To explain the construction, let  $n:N$ and $c:C$ be considered as input parameter (not as concrete objects), i.e. as pair 
$(n,c)$ of type $N\times C$.  We are going to follow the data-flow in the graph of $op$. 
\begin{itemize}
\item Operation $proj_{N,C}$ applied to $(n,c)$ returns two objects: \\ 
$proj^N_{N,C}(n,c)$ denoted by $n^0$ and 
$proj^C_{N,C}(n,c)$ denoted by $c^0$.
\item $Copy_N$ applied to $n^0$  returns two objects:
$Copy^1_N(n^0)$ denoted by $n^1$ and $Copy^2_N(n^0)$ denoted by $n^2$. 
\item $Copy_C$ applied to $c^0$ returns:  
$Copy^1_C(c^0)$ denoted by  $c^1$ and  $Copy^2_C(c^0)$, that is used again for copying.  
\item $Copy^1_C(Copy^2_C(c^0))$ is denoted by $c^2$, and   
$Copy^2_C(Copy^2_C(c^0))$ is denoted by $c^3$. 
\item $c^1, c^2, c^3$ are copies of $c^0$, and  $n^1, n^2$ are copies of $n^0$. Actually, they are copies of $c$ and $n$ respectively. 
\item Apply  $Suc$ to $n^2$, i.e.  $Suc(n^2)$. Then, copy the result twice, i.e.  $Copy_N$ applied to $Suc(n^2)$ returns:
$Copy^1_N(Suc(n^2))$ denoted by $\underline{n}^1$.    
$Copy^1_N(Copy^2_N(\underline{n}^1))$ is denoted by $\underline{n}^2$, and \\   
$Copy^2_N(Copy^2_N(\underline{n}^1))$ is denoted by $\underline{n}^3$. 
\item $\underline{n}^1$, $\underline{n}^2$ and  $\underline{n}^3$ are three copies of $Suc(n^2)$, i.e. they are the same as $n+1$. 
\item Let $apply_{((C;N) \rightarrow (A\rightarrow A)),(C;N)}(c^1; n^1)$ be denoted by $c^{1}_{n^1}$; it is the $n$-th element of sequence $c$, i.e. $c(n)$. 
\item Let $apply_{((C;N) \rightarrow (A\rightarrow A)),(C;N)}(c^2; \underline{n}^1)$ be denoted by $c^2_{\underline{n}^1}$; it is the same as $c(n+1)$. 
\item $compose_{A,A,A}(c^{1}_{n^1}; c^2_{\underline{n}^1})$ is  the composition of $c(n)$ and $c(n+1)$. Let it be denoted by $f$. 
\item Change in the sequence $c^3$,  the element $\underline{n}^2$-th  to $f$, i.e.  $Change_C(\underline{n}^2; f; c^3)$, and denote it by $\underline{c}$.  In this way $(n+1)$-th element element of sequence $c$ was changed to   $f$. 
\item Join  $\underline{n}^3$ and  $\underline{c}$ into the pair, i.e. $join_{N,C}(\underline{n}^3; \underline{c})$, and denote it by $(n+1, \underline{c})$.
\item Starting with  $(n,c)$ as the input in the construction we get  $(n+1, \underline{c})$ as the output. Actually, the only change that was made in the original sequence $c$ was to replace the $(n+1)$-th element of $c$ by the composition of two operations $c(n)$ and $c(n+1)$.
\end{itemize}
The description of the construction of the operation $op: N\times C \rightarrow N\times C$, in  Fig. \ref{rys6-7}, is completed. 

Let $N\times C$ be denoted by $D$. Then, $op$ is of type $D\rightarrow D$. 
Now, operation $Iter_D$ can be applied to $op$. 

Note that $Iter_D (n; op)$ is the $n$-th iteration of operation  $op$ such that $Iter_D (n; op) (1,c)$ is $(n+1, \underline{c})$ such that 
for any $k= 1,\ ...\  n+1$, the element $\underline{c}(k)$ is the composition of the first $k$ elements  of $c$. For $m$ greater than $n+1$,   element $\underline{c}(m)$ is the same as  $c(m)$. 
Note that $n:N$ is the parameter of the construction.  
To construct  Grzegorczyk's iterator, the operation $Rec_A$, shown as the second graph in Fig.  \ref{rys6-7}, must be constructed first. 
  
Let us follow the data flow in the graph of $Rec_A$ for the parameters $n:N$ and $c:C$.  
\begin{itemize}
\item At the node $Iter_D$, the iteration is applied only to the operation $op$ (constructed in Fig. \ref{rys6-7}) leaving the input  $n:N$ open, i.e. 
$Iter_D(*; op)$ is operation of type $N \rightarrow (D\rightarrow D)$. 
Note that this is amorphous partial application. The plug of this operation is of type $D\rightarrow D$, so that for parameter $n$, the operation $Iter_D(*; op)$ results in an operation (denoted by $g$) of type $D\rightarrow D$. \\
Note that here the unfolding is done in the graph for a concrete value of the parameter $n:N$, that is, the operation $op$ is iterated $n$-times.  
\item $join_{N,C}(1;c)$ is denoted by $(1,c)$. Note that $c$ is a parameter. 
\item $apply_{(D\rightarrow D), D}(g; (1,c))$ is denoted by $(\underline{n}, \underline{c})$. 
Let $proj^N_{N,C}(\underline{n}, \underline{c})$ be denoted by $\underline{n}$, and $proj^C_{N,C}(\underline{n}, \underline{c})$ by $\underline{c}$. 
\item $\underline{n}$ is the same as $n+1$, and for any $k=1,2,\ ...\ ,n+1$;  \ \  $\underline{c}(k)$ is the composition of the $k$ first operations in the original sequence $c$. 
\item $Pred(\underline{n})$ is the same as $n$. 
Finally, $apply_{C,N}(\underline{c}, Pred(\underline{n}))$ is operation of type $A\rightarrow A$. It is the composition of the $n$ first elements (operations) of the original sequence $c$.
\end{itemize}
This completes a description of the construction of the operation $Rec_A: (N;C)\rightarrow (A\rightarrow A)$  in Fig. \ref{rys6-7}. 

For the  inputs $n$ and $c$, the output, i.e. $Rec_A(n;c)$, is  the composition of the first $n$ elements (operations) from the input sequence $c$. This is the exact meaning of the construction of $Rec_A$. 

However, applying currying, the operation $Rec_A$ may be presented equivalently  as the operation  $\bar{Rec}_A$ of type $C\rightarrow (N \rightarrow (A\rightarrow A))$, i.e. of type $C\rightarrow C$. This may suggest that the operation $\bar{Rec}_A$ takes as its input a complete infinite sequence, and returns as the output also a complete infinite sequence. It is not true. By the construction of $Rec_A$ it is clear that $n:N$ is the parameter for this construction, i.e.  for any $n:N$, the construction is a finite structure referring only to the first $n$ elements of the input sequence $c$. 

Operation $\bar{Rec}_A$ corresponds to the Grzegorczyk's iterator. 
As an object, it can be used in more sophisticated constructions.

%%%%%%%
%%%%%%%%%%%
\subsection{ Summary of the level zero } 

The primitive type of natural numbers, the primitive constructors, and primitive operations are the basis for construction of objects that together   constitute the level zero of the Universe. 
The level may be viewed as a grounding (concrete semantics) for Grzegorczyk System, i.e. the Grzegorczyk's idea of {\em primitive recursive objects of all finite types}.  However, the general recursive objects (according to G\"{o}del-Herbrand definition) have grounding (as constructions) on higher levels.  Also a grounding for Girard System F needs higher levels of the Universe.   

The hardware interpretation of types, their constructors, and operations as (dynamic) configurations of connections between  plugs and sockets in nested boards, is important. This gives rise to comprehend the notion of type, the notion of object and its construction as a parameterized inductive structure that is finite  if the values of all parameters are fixed.  The  interpretation is in opposition to formal theories, where object is described as a term, construction amounts to substitution and lambda abstraction whereas computation to beta reduction.  

%%%%%%%%%%%%% 
%%%%%%%%%%%%%%%%%%%%%%%
\section{ Level 1 }
\label{level1}

Passing from level 0 to level 1 consists in handling types as objects. So that operations can be performed on types. The type constructors $+$, $\times$, 
%$||$, 
and $\rightarrow$ may be seen as operations. Also all primitive and complex operations indexed by types can be seen as operations level at 1, taking types as parameters. 

For example, primitive recursion schema (Grzegorczyk's iterator) $\bar{Rec}_A$ indexed by type $A$ may be abstracted to operation such that for arbitrary type $A$, the operation returns $\bar{Rec}_A$. It looks somehow as Lambda abstraction from System F. However, it is not a term. It is a concrete operation. 

Let $Types^0$ denote the type of all simple types constructed from the type of natural numbers by the constructors  $+$, $\times$, and $\rightarrow$. Then, the operation $\phi$ such that for any type $A$, \ \  $\phi (A)$ is $\bar{Rec}_A$ can be constructed, and has hardware interpretation. How can it be done? What is the type of this operation?  

As a primitive type, $Types^0$ is well defined with the primitive object $N$, and primitive operations $+$, $\times$, and $\rightarrow$. 
Actually, $Types^0$ is an inductive type. So that an operation can be constructed (on level 1) for enumerating all types from level 0. Let us fix one of such operations, and denote it by  $Ind^1: N \rightarrow  Types^0$. Hardware interpretation of this operation is that for a number $n:N$, it produces $Ind^1(n)$ as an appropriate nested board of sockets and plugs. 

Since the level 1 is the extension of level 0, type constructors, the primitive type $N$, primitive operations, and all types and objects that can be constructed at level 0, belong to level 1. 
In the same manner, the consecutive levels are introduced, so that at level 2 the type of all types at level 1 (denoted by $Types^1$) is introduced as a primitive type, and in general at level n+1, the type $Types^n$ is introduced as a primitive type.  

Introducing $Types^0$, as a primitive type at level 1, makes essential  change to its structure. The type constructors, as well as all objects and operations indexed by simple types become operations at level 1 that take simple types as their parameters. Well known new type constructors (see Martin-L\"{o}f Type Theory, System F, and CoIC) for dependent types emerge in the natural way.  

Let us start with type constructors $\times$, $+$, and $\rightarrow$. At the level 1, they can be considered as operations that take two types from $Types^0$ and return a complex type.  The operations may be presented as    
$
\times^1 
$,   
$
+^1 
$,
% $||^1 $,
 and 
$
\rightarrow^1 
$.
All of them are  
of type $(Types^0; Types^0 ) \rightarrow_1  Types^0$. 

Here  $\rightarrow_1$ is regarded as constructor at level 1, and is the extension of $\rightarrow$ from level 0, whereas $\rightarrow^1$ is an operation at level 1. Analogously, for the rest of constructors.  
Does it make sense? Although, it is formally correct, it is an unnecessary complication. Applying the well known parametricity (Reynolds  \cite{Reynolds}), the type constructors become generic primitive operations that for {\em any} two types produce a new complex type. Hence, the constructors are operations defined for all types at level 0, at level 1, at level 2, and so on, i.e. for all types at all levels. 
The parametricity may be also applied to any object indexed by types, so that it may be abstracted to generic operations defined on all types at all levels.   

For example,  $plus^1$ is an operation at level 1, i.e.  its two sockets are of type $Types^0$, whereas the type of its plug is determined by input objects. That is, for two input objects $A$ and $B$ (simple types from level 0) the result $plus^1(A; B)$ is the operation  $plus_{A,B}: (A;B) \rightarrow  (A + B)$. 
% Types of these operations are well known dependent types introduced in the next section. 
However, by applying the parametricity, $plus$ is a generic primitive operation defined (like the generic type constructors) for all types at all levels.  

It seem that it is reasonable to consider the ``the super type'' of all types denoted by $Types$. Then, the generic type constructors ($\times$, $+$, and $\rightarrow$) are of type $(Types; Types)\rightarrow Types$. 
Does it make sense? Is $(Types; Types) \rightarrow Types$ an object of the ``the super type'' $Types$? For the formal theories (like Martin-L\"{o}f Type Theory and System F) introducing ``the super type'' results in terms that have no normal form, so that from the computational point of view they have no meaning. However, the approach presented in the paper is not yet another formal theory. It is to be grounded in hardware. 

The super type $Types$ as a completed object has no grounding. It is only a useful abstract notion. The general principle of the proposed approach is that the Universe is never completed, and in fact at any level of its construction it is a finite structure that has a hardware interpretation. Construction of the Universe is a never ending process, so that at any moment of time of the process, only a finite number of types and finite number of objects are constructed. Moreover, for any of such types and any of such objects, if its construction is inductive, then (for that moment of time) the construction is only partial up to fixed values of the inductive parameters. 
At any moment of time of the construction process of the Universe, the current level is finite, say k, so that the super type $Types$ is interpreted as $Type^k$. For this interpretation of $Types$, it is only a useful abstract notion, and there is no contradiction.  
  
Hence, the following type constructors and primitive operations can be interpreted as generic and defined on the super-type $Types$, i.e. defined for all types at all levels. 
\\   
$+$, $\times$, $\rightarrow$, 
% $||$,  
$join$,
$proj$, 
$plus$,
$get$,
$compose$,
$apply$,
$const$,
$id$, 
$Copy$, 
$Iter$,
$Change$,  $currying$, $uncurrying$

For any of the above operations, the type of output object is determined by the input objects. This very determination is given by an operation $F: D \rightarrow  Types$ different for different operations. For operation $plus$,  the socket of $F$ is $(Types;  Types)$ so that $F(A;B)$ is $A+B$. In this way, $plus(A;B)$ is the same as $plus_{A,B}: (A;B) \rightarrow  (A + B)$. 

Types of the above operations are well known dependent types, Girard \cite{Girard71} and Martin-L\"{o}f  \cite{Martin-Lof}. There are two new type constructors:  $\Pi$ as the generalization of arrow $\rightarrow$, and  $\Sigma$ as the generalization of disjoin union $+$.  

%%%%%%% 
\subsection{ Dependent types }
\label{typy-zalezne}

For operation $F: A \rightarrow  Types^0$, an object of type  $\Pi_A F$ is the operation $g$, such that for any  $a:A$,\ \  $g(a)$ is of type $F(a)$. 
If $F$ is a constant, i.e. its value is $B$, then the type $\Pi_A F$ is the same as $A \rightarrow B$. 

Object of type $\Sigma_A F$ is of the form $(a; b)$, where $a:A$ and $b: F(a)$. 

Operation $F$ may have multiple input (several sockets). 

Dependent type constructors are generic operations. So that their generic forms are denoted by $\Pi$ and $\Sigma$. 
For simplicity, $\Sigma_A F$ is denoted by $\Sigma F$, and $\Pi_A F$ is  denoted by $\Pi F$, if from the context it is clear that $F$ is of type $A \rightarrow Types$.

%%%%%%%%%%%%%%% IMPORTANT <<<<<
% Their types are determined by the  generic and primitive type constructor $\psi$ of type $Types \rightarrow Types$ such that for any $A$ , \ \ $\psi(A)$ is  $A\rightarrow Types$.  
% Then, $\Pi$ is of type $(Types; \Pi_{Types} \psi) \rightarrow Types$ such that  for any $A:Types$ and any $\mathbb{F}: \Pi_{Types} \psi$,  the type $\Pi(A; \mathbb{F})$ is  $\Pi_A \mathbb{F}(A)$.  
%
% Analogously, $\Sigma$ is of type $(Types; \Pi_{Types} \psi) \rightarrow Types$ such that for any  $A$ and any $\mathbb{F}$, the type $\Sigma(A; \mathbb{F})$ is  $\Sigma_A \mathbb{F}(A)$.  
%
% Note that in the definitions of the type of $\Pi$ and $\Sigma$, the constructor $\Pi_{Types}$ was used. However, there is no circularity if the type is interpreted in finite structures in the same way as the interpretation of type $(Types; Types) \rightarrow Types$.  
%
%For simplicity,  $\Sigma(A; \mathbb{F})$ will be denoted by $\Sigma F$, and $\Pi(A, \mathbb{F})$ will be denoted by $\Pi F$, if from the context it is clear that $F$ is $\mathbb{F}(A)$. 
%
% >>>>>> 
%
%For the interpretation of $\Pi$ and  $\Sigma$ as logical quantifiers (see the next section) the first argument (type $A$) indicates the scope of the quantification. 

{\bf Constructor of objects of type $\Sigma F$ from objects of type $\Pi F$} . If for $a:A$ an object $b$ of type $F(a)$ is constructed, then, in fact, also an object of type $\Sigma F$ is constructed. 
For operation $F:A \rightarrow  Types$  and operation  $f: \Pi F$ (i.e. for all $a:A$, \ \ $f(a): F(a)$) new constructor $\sigma_F$ of type $\Pi F \rightarrow  (A \rightarrow  \Sigma F)$ is introduced.  Operation  $\sigma_F(f): A \rightarrow  \Sigma F$ is such that for any  $a:A$, \ \ $\sigma_F(f)(a)$ is of type $\Sigma F$, i.e. it is of the form $(a; f(a))$. The new constructor  also has its generic form denoted by  $\sigma$. 

Hardware interpretation of the dependent types is no so easy. The operation $F$ must be taken as the basis. The type $\Sigma F$ is interpreted as a board consisting of type $A$ and operation $F$ that dynamically determines (for any $a:A$) the type $F(a)$ in the board.   For an object of the type $\Pi F$, say $f$, its socket is of type $A$ whereas type of the plug is dynamically determined. So that, the type $\Pi F$ may be interpreted as a board consisting of the socket of type $A$, and operation $F$ that dynamically determines the type of the plug. 

It seems that the above interpretation make sense. However, any new type constructors should be accompanied by object constructors for these new types.  Operation $F$ alone does not indicate how to construct objects of type $\Pi F$, and objects of type $\Sigma F$. However, the objects of dependent types emerged in the natural way during the construction of the Universe, as generic type constructors and generic parametric versions of the primitive operations. 

Summing up the level 1 of the Universe, generic primitive constructors and operations have been introduced that are defined for all types at all levels. Although their grounding is (and must be) always in a finite hardware structure being the current state of construction process of the Universe, the constructors and the operations are useful abstract notions. Let us recall that the main purpose of the paper is a design of generic mechanisms for managements of connections in large arrays of hardware functional units (first order functions). From this point of view, the abstractions are very useful tools.  

However, these abstractions do not exhaust the toolkit. There are also relations, and important generic operations corresponding to {\tt if-then-else} statement in programming.  

%%%% 
%%%%%%%%
\section{ Relations }
\label{relacje} 

Usually, binary relation is defined as a collection of ordered pairs of objects. This set theoretical definition is not sufficient.  Relation is (like  operation) a primitive notion. It seems that it corresponds to a primeval generic method of comparing two objects. 
For any primitive type there is at least one elementary binary relation between objects of this type.

Well known equality types in Martin-L\"{o}f Type Theory \cite{Martin-Lof}, when parametrized for a fixed type, give an equality relation on that type. However, this relation is pure syntactical one, and their evaluation is based on term reductions to canonical normal forms. 

%%%%%%
\subsection{ Primitive relations on natural numbers  }   
\label{el-rel-nat}

The relations $Equal_N$, $Lesser_N$, and $Greater_N$ are primitive relations on $N$ with the following grounding (hardware interpretation) presented below.  

There are two sockets of type $N$, one denoted by $N'$ for $n$ , and one denoted by $N''$ for $k$. It is supposed that for each of the sockets its state can be evaluated as either empty or not empty. This evaluation may be considered as the most primitive relation (property) for the type of natural numbers. 

Put the signal $n:N$ into the socket $N'$, and the signal $k:N$ into the socket $N''$.\\
{\bf Procedure N.}  Check the states the two sockets. Until one of the sockets is empty, apply to each of them $Pred$, that is, $Pred(n)$, and $Pred(k)$. It means to remove from each of the sockets one elementary signal; then, go to the beginning of the procedure. If both sockets are empty, then this is the witness (proof) for the proposition $Equal_N(n;k)$ to be true. If the socket $N'$ is not empty and the socket $N''$ is empty, then it is the witness (proof) for the proposition $Greater_N(n;k)$ to be true. If the socket $N'$ is empty and the socket $N''$ is not empty, then it is the witness (proof) for the proposition $Lesser_N(n;k)$ to be true. 

For any $n:N$ and $k:N$, 
$Equal_N(n;k)$ is a primitive proposition corresponding to the states of the two sockets with the intuitive meaning that $n$ is equal to $k$. 
Analogously for $Lesser_N(n;k)$ with the meaning that $n$ is lesser than $k$, and for $Greater_N(n;k)$ with the meaning that $n$ is greater than $k$. 

Note that the results of evaluations in the Procedure N are witnesses (proofs) that correspond the intuitive notion of {\em truth}. That is, a proposition is {\em true} if there is a corresponding witness for this proposition. If there is no witness, then the proposition is {\em false}. 

The relations seem to be operations. If it is so, then what are their plugs? 
According to the famous idea of {\em Curry-Howard  propositions-as-types}, the proofs are objects of propositions considered as types.  
Note that the grounding of the propositions is in the Procedure N. So that, a single proposition can not be considered separately without reference to the rest of the propositions of the Procedure N.  Some of the propositions are false, so that the corresponding types are not inhabited, i.e. are empty. So far all introduced types were inhabited. 

Note that the type constructors ($+$, $\times$, $\rightarrow$, $\Pi$ and  $\Sigma$) have also logical interpretation as disjunction, conjunction, implication, and quantifiers. However, the emptiness of a proposition, considered as a type, causes severe problems. One of them is the hardware interpretation of empty type. Empty type is nonsense, it can not be realized, it does not exist. If $A$ is an empty type, then also the type $A \rightarrow B$ is empty contrary to the classical logical interpretation of $A \rightarrow B$ as implication that must be true in this case. 

Although the notion of type and the notion of proposition have a lot of common, they should be separated. Let $Prop$ denote the sort of propositions.  The same distinction was made in CoIC \cite{Coq}, however there the sort $Prop$ corresponds merely to a formal logic and its formulas and terms. 
The hardware interpretation of $Prop$ requires an introduction of primitive propositions, and a procedure for their evaluation. For natural numbers these primitive propositions (and their hardware interpretation) are introduced by the Procedure N. 

Hence, for any type $A$, the ``type'' $A \rightarrow Prop$ has sense, so that  operations of this ``type'' may be constructed. They are called relations. It seems that this ``relation type'' is special and the usual type constructors can not be applied to it. 
Any of such relations (say $R:A \rightarrow Prop$) can be quantified by applying $\Pi$ and $\Sigma$ to this relation, so that $\Pi R$ and $\Sigma R$ are propositions of the sort $Prop$.  
General form of relation type is  $(A_1; ... ; A_k) \rightarrow Prop$.
Complex relations can be constructed from the primitive relations by using constructors  $+$, $\times$, $\Pi$, $\Sigma$ (interpreted in $Prop$), and negation introduced in the next section \ref{kan}.  Although the constructors are interpreted in $Prop$, following the {\em Curry-Howard  propositions-as-types}, and for notational simplicity, they are denoted by the same symbols.

%%%%%%%%
\subsection{ Complex relations and complex propositions }
\label{kan}

\begin{figure}[h]
	\centering
	\includegraphics[width=0.6\textwidth]{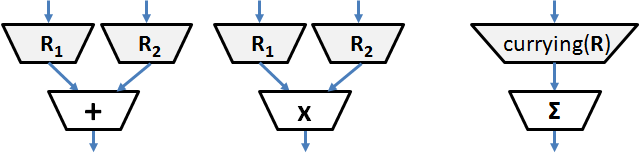}
	\caption{ Disjunction and conjunction of $R_1$ and $R_2$, and composition of $currying(R)$ of type $A\rightarrow (B \rightarrow Prop)$  with $\Sigma$ }
	\label{currying(R)}
\end{figure}

For two relations  $R_1:A \rightarrow Prop$, and  $R_2:B \rightarrow Prop$, their conjunction and disjunction is constructed in the following way, also shown in Fig. \ref{currying(R)}.  
Note that in the logical interpretation, the constructors $+$ and $\times$ are of type $(Prop; Prop) \rightarrow Prop$.  

Disjunction is denoted by $R_1 + R_2$, and is amorphous double composition of $R_1$, $R_2$, and operation $+$. The composition is done by connecting the plug of one relation to the one of the sockets of $+$, and by connecting the plug of the second relation to the second socket of $+$. 

Conjunction is denoted by $R_1 \times R_2$, and is constructed in the very similar way.  

$R_1 \times R_2$ and $R_1 + R_2$ are relations of type $(A; B) \rightarrow Prop$. 

%dodatek WAŻNY 17 marca 2017 -- sprawdź 
For two relations  $R_1$, and  $R_2$, both of type $A \rightarrow Prop$, the special generic forms of conjunction and disjunction (corresponding to informal $(R_1(a) + R_2(a))$  and $(R_1(a) \times R_2(a))$) are constructed in the following way. 

Operation $Copy_A: A \rightarrow (A;A)$ with its two plugs of type $A$ is composed amorphously with $R_1$ and  $R_2$ by connecting the plugs to the sockets of  $R_1$ and  $R_2$. Then the resulting operation of type $A \rightarrow (Prop; Prop)$ is composed again amorphously with the conjunction operation $\times : (Prop; Prop) \rightarrow Prop$ (resp. disjunction $+ : (Prop; Prop) \rightarrow Prop$). The final operation is denoted (a bit informally) by $[R_1 \times R_2]$ (resp. $[R_1 + R_2]$) both of type  $A \rightarrow Prop$ such that for any $a:A$,  \ \ $(R_1(a) + R_2(a))$  and $(R_1(a) \times R_2(a))$ are generic forms of $[R_1 \times R_2](a)$ and $[R_1 + R_2](a)$.
 
Negation of a single separate proposition does not have sense. For example, $\neg Equal_N(1;2)$ (where $\neg$ is negation constructor) makes sense (see the Procedure N) only in the presence of $Greater_N(1;2)$ and $Lesser_N(1;2)$. Hence, the negation for relation can be grounded only if there are complementary relations. For $Equal_N$, the complementary relations are $Greater_N$ and $Lesser_N$. Hence, the negation of a relation is its complement.  

$\neg Equal_N$ (i.e. the complement of $Equal_N$) is $[Greater_N + Lesser_N]$. 
 
 $\neg Greater_N$ (i.e. the complement of $Greater_N$) is $[Equal_N + Lesser_N]$. 
 
 $\neg Lesser_N$ (i.e. the complement of $Lesser_N$) is $[Equal_N + Greater_N]$. 
 
Hence, the constructor $\neg$ can be applied only to a relation whose complement has already been constructed. This very complement is the negation of the relation.  Since negation is complement, the meaning of  
 double negation, i.e. $\neg\neg R$ is $R$. 
% ??? And for two relations $R_1$ and $R_2$ (having complements, i.e. resp. $\neg R_1$ and $\neg R_2$) the De Morgan's laws hold, that is, $\neg(R_1 + R_2)$ is defined as  $(\neg R_1 \times \neg R_2)$, and  $\neg(R_1 \times R_2)$ is  defined as $(\neg R_1 + \neg R_2)$.  

%

Note that the quantifiers $\Pi_A$ and $\Sigma_A$ are of type $(A \rightarrow  Prop) \rightarrow Prop$, i.e. for any relation  $R: A \rightarrow  Prop$, \ \ $\Pi_A R$ and $\Sigma_A R$ are propositions belonging to the sort $Prop$. Bottom index $A$ will be frequently omitted if from the context it is clear what it is.  

For any relation $R$ of type $A \rightarrow  Prop$, the constructors $\Pi$ and $\Sigma$ may be applied resulting in two propositions $\Pi R$ and  $\Sigma R$. The first proposition corresponds to the formula $\forall_{x:\bar{A}} \bar{R}(x)$ in formal logic, the second proposition corresponds to formula $\exists_{x:\bar{A}} \bar{R}(x)$, where $\bar{A}$ and $\bar{R}$ are symbols in formal language of the logic, corresponding to type $A$ and relation $R$, and $x$ is a variable. 

If a relation has multiple input (several sockets) like  $R:(A;B) \rightarrow Prop$, then $\Pi R$ corresponds to formula $\forall_{x:\bar{A}} \forall_{y:\bar{B}} \bar{R}(x;y)$; analogously for  $\Sigma$. 

The formula $\forall_{x:\bar{A}} \exists_{y:\bar{B}} \bar{R}(x;y)$ corresponds to the proposition that is constructed in the following way. 
First, relation $R$ must be  transformed by $currying$ to the operation  $currying(R): A\rightarrow (B \rightarrow Prop)$, such that $R(a;b)$ is the same as $currying(R)(a)(b)$. Note that $\Sigma_B$ is of type $(B \rightarrow Prop) \rightarrow Prop$. 
The amorphous composition of $currying(R)$ and $\Sigma_B$, see Fig. \ref{currying(R)}, results in the operation of type $A \rightarrow Prop$  denoted by $S$. 
So that  $\Pi_A S$ corresponds to the formula $\forall_{x:\bar{A}}\exists_{y:\bar{B}} \bar{R}(x,y)$.
Actually, the construction described above is generic and can be generalized to any types and any relations.   

% <<<<<<<< to co niżej można pominąc w skroconej wersji
For relation $R: A \rightarrow Prop$, application of $Copy$ to $R$ (i.e.  $Copy(R)$) returns two the same operations  $Copy^1(R)$ and $Copy^2(R)$. Suppose that the negation of $R$ (its complement) has been already constructed.
 
Relation $Copy^1(R) + \neg Copy^2(R)$ is of type $(A;A') \rightarrow Prop$. Amorphous composition of operation $Copy_A$ (of type  $A \rightarrow (A;A')$ ) with  $Copy^1(R) + \neg Copy^2(R)$  gives the relation of type $A \rightarrow Prop$ that corresponds to {\em tertium non datur} (TND) in the formal logic. Denote this relation by $R^{TND}$. For any $a:A$, \ \ $R^{TND}(a)$ is, informally, $R(a) + \neg R(a)$. 
 Since $\neg R$ is the complement of $R$, then the proposition $R^{TND}(a)$ is true for any $a:A$. Note that this holds only for a relation having already constructed complement.    
% >>>>>>>>>>>

In order to summarize the Section devoted to relation, let us come back to the Procedure N, and primitive relations $Equal_N$, $Lesser_N$ and $Greater_N$. It determines important properties concerning the primitive relations, primitive operations $Suc$ and $Pred$. These properties (as true propositions) should be considered as axioms for the type $N$. For example, if $Equal_N(n;k)$ is true, then also  $Lesser_N(n;Suc(k))$ and $Greater_N(Suc(n);k)$ are true. %, see section \ref{aksjomaty}. 

Also the following informal proposition \\ 
{\em for all $n$ and for all  $k$: either $Equal_N(n;k)$ or $Lesser_N(n;k)$ or $Greater_N(n;k)$} 
\\ 
is true. Sentences belonging to the sort $Prop$ can be constructed on the basis of the above constructions (presented in this Section) of complex relations and propositions. As the axioms, characterizing completely the Procedure N, they have proofs that are to be introduced as primitive objects along the primitive relations corresponding to the type $N$ of natural numbers. The proofs can be used to constructions of another proofs of complex sentences. These very axioms correspond to the axioms of the formal theory of Arithmetics, i.e. Peano Arithmetics.

%%%%%%%%%
%%%%%%%%%%%%%%%
\section{Conditions  }

The notion of condition and its verification is introduced. 
\begin{figure}[h]
	\centering
	\includegraphics[width=0.95\textwidth]{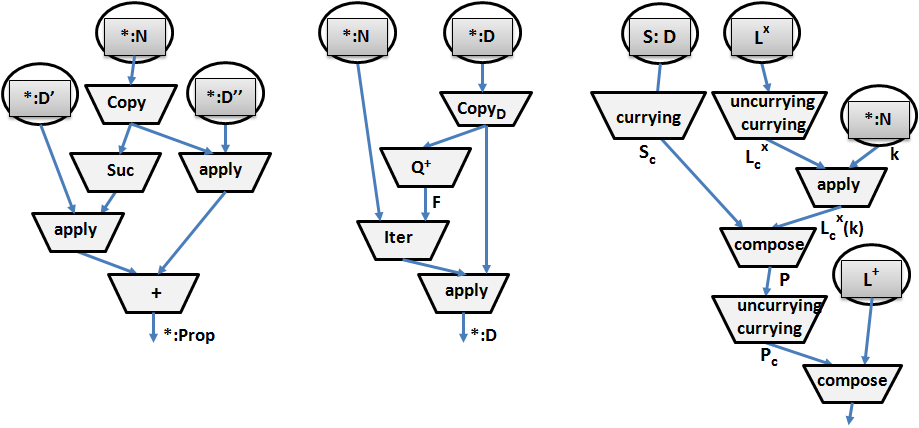}
	\caption{On the left, auxiliary operation $op$, and operation $L^+$. On the right, complex condition where  $E$ denotes type $(N'; N)\rightarrow Prop$ }
	\label{rys17}
\end{figure} 
Generally, a condition consists of disjunctions and conjunctions of primitive propositions and their negations.  The disjunctive normal form (disjunction of conjunctions) is very convenient for verification, i.e. once one component of the disjunction is verified as true, then the condition is true; if one component of a conjunction is false then the conjunction is false.  A generic  verification method can be constructed on the basis of the primitive propositions. 
Simple examples of informal parametrized conditions for relation $R:N \rightarrow Prop$ are as follows. 
\begin{itemize}
\item {\em exists $i$ such that $k \leq i \leq k+n$ and $R(i)$} 
\item {\em for all $j$, if $l \leq i \leq l+m$, then $R(j)$}
\end{itemize}
Note that here the phrases {\em ``for all''} and {\em ``exists''} do not correspond to the constructors $\Pi$ and $\Sigma$. 

The parametrized conditions (conditions, for short) correspond to operations of type 
$(N; N; (N \rightarrow Prop)) \rightarrow Prop$. 
Disjunction $+$, conjunction $\times$, and negation $\neg$ are needed to construct the conditions. 
 
Let $D$ denote $N \rightarrow Prop$. Type $D'$ denotes a copy of $D$. We are going to construct the operation that corresponds to the first of the above conditions. 
Auxiliary operation $op$ of type $(D'; N; D'') \rightarrow Prop$  is constructed in Fig. \ref{rys17}. For $R_1$, $R_2$ and $i:N$,  \ \ 
$op(R_1; i; R_2)$ is the same as $R_1(i+1) + R_2(i)$. 
Applying $currying$ and $uncurrying$ to operation $op$ we get operation $Q^+: D'' \rightarrow (D'\rightarrow (N \rightarrow Prop))$ such that  
$Q^+(R_2)(R_1) (i)$ is the same as $R_1(i+1) + R_2(i)$. 

We are going to construct operation  $L^+$ of type $(D;N') \rightarrow (N \rightarrow Prop)$ such that for any $R:D$, $n:N'$, and $k:N$, \ \  $L^+(R;n)(k)$
 corresponds to 
$(R(k) + R(k + 1) + R(k + 2 ) + ... +  R(k+n))$, i.e. informally  {\em (exists $i$ such that $k \leq i \leq (k+n)$ and $R(i)$)}. 
The construction of operation  $L^+$ is shown in Fig. \ref{rys17}. 

Relation $R:D$ is copied, and  $Q^+$ is applied to $Copy^1(R)$. The result  denoted by $F$ is of type $D'\rightarrow D$. So that iteration $Iter_D$ can be applied to $F$, and $Iter_D(n; F)$ is of type $D' \rightarrow D$. Note that $F(R_1)(k)$ is $R_1(k+1) + R(k)$, and  $(F\circ F)(R_1)(k)$ is $R_1(k+2) + R(k+1) + R(k)$, whereas $(F\circ F\circ F)(R_1)(k)$ is $R_1(k+3) + R(k+2) + R(k+1) + R(k)$. So that  $Iter_D(n; F)(R_1)(k)$ is $R_1(k+n) + R(k+n-1) + ... + R(k+1) + R(k)$.  
Finally, applying $Iter_D(n; F)$ to $Copy^2(R)$, we get the required relation $L^+(R;n)$ of type $N \rightarrow Prop$. 

If in the construction of $L^+$ (more precisely in $Q^+$), the  operation $+$ is changed to $\times$, then the resulting  operation is denoted by $L^{\times}$.  
Then  $L^{\times}(R; m)(l)$ corresponds to $(R(l) \times R(l+1) \times ... \times  R(l+m))$, i.e. informally {\em (for all $j$, if $l \leq i \leq l+m$, then $R(j)$)}. 

Now, let $R$ be a relation having two variables, i.e. of type $(N'; N)\rightarrow Prop$. The complex condition (where $l, m$ and $k, n$ are parameters) corresponding to  \\  
{\em (there is $i$ and $l \leq l \leq l+m$ such that for all $j$ if  $k \leq j \leq k+n$ then $R(i;j)$)} \\ 
is constructed below, see also Fig. \ref{rys17} where in the rightmost  graph, $uncurrying$ and $currying$ are aggregated into one operation. 
  
Operation $L^{\times}$ is of type  $((N\rightarrow Prop); N'')\rightarrow (N'''\rightarrow Prop)$, where, argument $j$ corresponds to  $N$, argument $m$ corresponds to $N''$, and argument $l$ corresponds to $N'''$.  

Relation $R$ is of type $(N'; N)\rightarrow Prop$, where socket $N'$ corresponds to argument $i$, and socket $N$ corrsponds to argument $j$. 
Applying $currying$ to $R$, we get $R_c$ typu $N' \rightarrow (N\rightarrow Prop)$.

The plug of operation $R_c$ is of the same type as the type of one of the sockets of operation $L^{\times}$. These two operation can be composed, i.e.  $compose(R_c; L^{\times})$. Let the composition be denoted by $P$; it is of type $(N'; N'')\rightarrow (N'''\rightarrow Prop)$. 

Proposition $P(i;m)(l)$ is the same as $(R(i;l) \times R(i; l+1) \times ... \times  R(i; l+m))$. 

Applying $uncurrying$ and $currying$ to $P$, we get operation  $P_c$ of type $(N''; N''')\rightarrow (N'\rightarrow Prop)$. 

For the operation 
$L^+:((N'\rightarrow Prop); N'''')\rightarrow (N'''''\rightarrow Prop)$, argument $i$ corresponds to $N'$, argument  $n$ corresponds to $N''''$, and argument $k$ corresponds to $N'''''$. 
 
The plug of the operation $P_c$ is of the same type as the type of one of the sockets of operation $L^+$. 
The composition (by connection the plug with the socket) is of type 
$(N''; N'''; N'''')\rightarrow (N'''''\rightarrow Prop)$, and is denoted by $F$. The proposition $F(m;l;n)(k)$ is the same as 
 \\ 
 $((R(k;l) \times R(k; l+1) \times ... \times  R(k; l+m)) + (R(k+1;l) \times R(k+1; l+1) \times ... \times  R(k+1; l+m)) + ... + (R(k+n;l) \times R(k+n; l+1) \times ... \times  R(k+n; l+m)))$. 

More complex conditions can be constructed on the basic of $L^{\times}$ and $L^+$.

%%%%%%%%
%%%%%%%%%%%%%%%%
\section{Example of pure function-level programming }
\label{trees}

One may ask how these higher order types, operations and relations   correspond to real computations and programming. 
They are mathematical objects with clear hardware grounding, just as John Backus \cite{Backus} postulated. 
The example presented below shows that programming on such mathematical objects is possible. It also shows that the higher order operations may be used in programming directly without terms (for denoting the operations) and rewriting rules. 

Let us recall that $const_{N,N}$ denotes the constant operation of type $N\rightarrow(N \rightarrow N)$ such that 
for any $n_c:N$, \ \ $const_{N,N}(n_c): N \rightarrow N$ is such that for  any $k:N$,  \ \ $const_{N,N}(n_c)(k)$ is $n_c$. 

The programming example is simple and is based on the three operations: ${\tt node}$, ${\tt father}$ and  ${\tt leaf}$ all of type $N \rightarrow N$, that along with the value of the parameter $n$ are interpreted as a data structure called tree. The tree structure may be modified by two operations: $add$ for adding new node to the tree, and $del$ to remove a leaf node from the tree. 

The parameter $n:N$ denotes the current scope of the construction of the data structure. 
\begin{itemize}
\item  ${\tt node}: N \rightarrow N$.  For any $i:N$, \ ${\tt node}(i)$ is either $1$ (denotes an already constructed  {\tt node}),  or $2$ (denotes  deleted {\tt node}), or $3$ (and greater) denotes unspecified node outside of the current scope of the construction. Initially, ${\tt node}$ is $Change_{N}(1;1;const_{N,N}(3))$, i.e. ${\tt node}(1)$ is set as $1$ (it is the root), and  for  $i$ grater than $1$, ${\tt node}(i)$ is $3$.

\item ${\tt father}: N \rightarrow N$. ${\tt father}(i)$ is interpreted as the node that is the father of  node $i$.  Initially, it is the constant operation $const_{N,N}(1)$.  
The {\tt node} and the parameter $n:N$ determine which inputs of ${\tt father}$ have intended meaning, i.e.  ${\tt father}(k)$ is meaningful (in the tree structure) for $k$ not greater than $n$ and if $k$ is not a removed node. 
\item ${\tt leaf} :N \rightarrow N$. 
 For any $i:N$, \ \ ${\tt leaf}(i)$ is either $1$ (it is a leaf if ${\tt node}(i)$ is also $1$), or $2$ (is not a leaf), or $3$ (and greater) as not constructed or already deleted. 
Initially, ${\tt leaf}$ is $Change_{N}(1;1;const_N(3))$.    
\end{itemize}
The initial operations: ${\tt node}$, ${\tt father}$ and ${\tt leaf}$, and the parameter $n=1$ represent the tree consisting only of the root. Applying simple operations $add$ and $del$ constructed below,  arbitrary complex tree structures can be constructed.  

The parameter $n:N$ denotes the number of the last constructed {\tt node}; initially it is set as $1$. The  next natural number  $Suc(n)$ is for the next {\tt node} to be constructed in the tree. 

In the construction of the operations, the parameter $n:N$  determines the current  scope  of the operations. For a number grater than $n$ (outside of the current scope), the nodes (that could have such numbers) are still not constructed, so that the reference to them does not have the intended meaning. 
Let $A$ denote $(N \rightarrow N)$, and $B$ denote  $(N; N'; A; A'; A'')$. 
Two operation are constructed to modify a tree; 
$add$ and $del$ both of type  $B \rightarrow B$. 
For input $(n; o; {\tt node}^{in}; {\tt father}^{in}; {\tt leaf}^{in})$ they return output $(k; o; {\tt node}_{out}; {\tt father}_{out}; {\tt leaf}_{out})$. 

Operation $add$ adds a node to the tree. The new node is given number $(n+1)$ and its father is an already existing node $o$.  
Operation $del$ removes node $o$ if it is a leaf. 

The following pseudo-code describes the operations. 

Operation $add$: 
\begin{enumerate}
\item  
{\bf if}  ($ Greater_N(o; n)) +  Equal_N({\tt node}^{in}(o);2)$) is true, i.e. $o:N$ is either outside of the current scope or it is a deleted node \\ 
{\bf then} do nothing; \\
{\bf else}  
%begin-else 1
	\begin{enumerate}
	\item Construct a new node $Suc(n)$ to be a child of the node $o$. 	 			 	That is, ${\tt node}(Suc(n))$ becomes $1$, i.e. 
	$Change_N(Suc(n); 1; {\tt node}^{in})$ 
	\item  ${\tt father}^{in}(Suc(n))$ becomes $o$, i.e.  			 			
	$Change_N(Suc(n); o; {\tt father}^{in})$ 	
	\item  ${\tt leaf}^{in}(Suc(n))$ becomes $1$, i.e. $Change_N(Suc(n); 1; {\tt 				leaf}^{in})$ denoted by ${\tt 	leaf}^{in}_1$

% composition 
\item {\bf if} $Equals_N({\tt leaf}^{in}_1(o);1)$, i.e. $o$ was a leaf in the tree \\
{\bf then} 
	\begin{itemize} %begin-then 2
	\item ${\tt leaf}^{in}_1(o)$ becomes $2$, i.e. $Change_N(o; 2; {\tt leaf}^{in}_1)$ 
	\end{itemize} %end-then 2
{\bf else} do nothing. 
\end{enumerate} %end-else 1
\end{enumerate}
Note that the phrase {\em 'do nothing'} corresponds to the operation  $id_B: B \rightarrow B$ such that $id_B(b)$ is $b$ for all $b:B$. 
 
The constructor {\bf if-then-else} must be a new primitive. It needs a condition, and two operations. Formally, it is introduced in the next subsection. %The input of the condition and the inputs of the two operations are the same.  

The first informal condition of $add$ is $( Greater_N(o; n)) +  Equal_N({\tt node}^{in}(o);2))$. However, $o$, $n$,  and ${\tt node}^{in}$ are parameters here. So that there is relation, denoted by $S_{11}$, of type $(N; N'; (N \rightarrow N)) \rightarrow Prop$ such that $S_{11}(o; n;  {\tt node}^{in})$ corresponds  to this informal condition. 

\begin{figure}[h]
	\centering
	\includegraphics[width=0.89\textwidth]{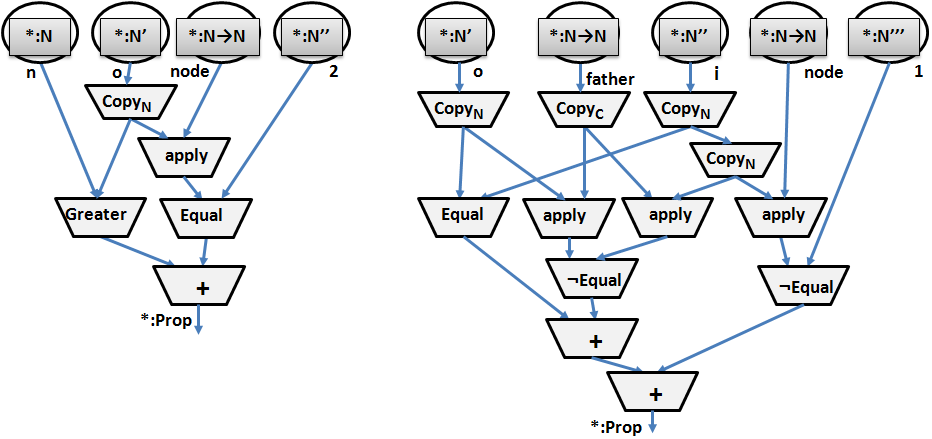}
	\caption{From the left, relations  $R_{11}$, and $R_{22}$}
	\label{rys11-22}
\end{figure}

To construct $S_{11}$, the relation $R_{11}$ is needed that is shown in Fig. \ref{rys11-22}. It is of type $ (N; N'; (N \rightarrow N); N'') \rightarrow Prop$ where the socket $N$ corresponds to $n$, the socket$N'$  to  $o$, the socket $(N \rightarrow N)$ to ${\tt node}^{in}$, and the socket $N''$ to $2$. 
\ \ The relation $R_{11}(*;*,*;2)$ is of type $(N; N'; (N \rightarrow N)) \rightarrow Prop$, and it is the required condition $S_{11}$. 

\begin{figure}[h]
	\centering
	\includegraphics[width=0.89\textwidth]{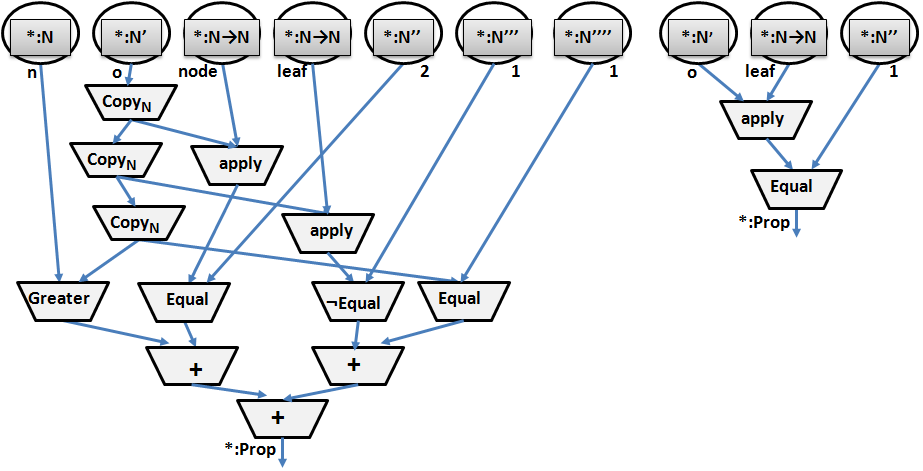}
	\caption{Relation $R_{21}$ and relation $R_{12}$}
	\label{rys21-12}
\end{figure}

The second informal condition of $add$ is $Equals_N({\tt leaf}^{in}(o);1)$. In Fig. \ref{rys21-12}, relation $R_{12}$ is constructed. It is  of type $(N';(N \rightarrow N); N'') \rightarrow Prop$, where socket $N'$ corresponds to $o$, socket $(N \rightarrow N)$ to ${\tt leaf}^{in}$, and socket $N''$ to $1$.  
\ \ The relation $R_{12}(*;*;1)$ of type $(N';(N \rightarrow N)) \rightarrow Prop$, is the required relation $S_{12}$ such that $S_{12}(o; {\tt leaf}^{in})$ corresponds to the second informal condition. 

~\\
Operation $del$: 
\begin{enumerate} 
 \item {\bf if} the condition $( Greater_N(o; n) +  Equal_N({\tt node}^{in}(o);2) + \neg Equal_N({\tt leaf}^{in}(o); 1) +  Equal_N(o;1))$,  is true,  i.e. either $o$ is outside of the scope, or $o$ is a removed  node, or $o$ is not a leaf, or $o$ is the root of the tree, \\
{\bf then} do nothing\\ 
{\bf else}   
	\begin{enumerate}	
	\item  ${\tt node}(o)$ becomes $2$, i.e. 
	$Change_N(o; 2; {\tt node}^{in})$
	\item ${\tt leaf}(o)$ becomes $3$, i.e. $Change_N(o; 3; {\tt leaf}^{in})$ is 			denoted by ${\tt leaf}^{in}_1$    
	
% composition 
\item {\bf if}  node $o$ is the only child of its  father, i.e. 
for all $i:N$ such that $i$ is not greater than $n$, either  $i$ is not a node, i.e. $\neg Equal_N(node^{in}(i); 1)$, or   
$(Equal_N(o; i) + \neg Equal_N({\tt father}^{in}(i); {\tt father}^{in}(o)))$, i.e. either the nodes $o$ and $i$ are the same, or they have different fathers  \\
{\bf then} 
	\begin{itemize}
	\item ${\tt leaf}^{in}_1({\tt father}^{in}(o))$ becomes $1$, i.e. 		
	$Change_N({\tt father}^{in}(o); 1; {\tt 	leaf}^{in}_1)$
	\end{itemize}
{\bf else}  do nothing   
\end{enumerate}
\end{enumerate}
The first informal condition in $del$ is 
 $( Greater_N(o; n) +  Equal_N({\tt node}^{in}(o);2) + \neg Equal_N({\tt leaf}^{in}(o); 1) +  Equal_N(o;1))$. In order to construct the relation (denoted by $S_{21}$) that corresponds to this condition, the auxiliary relation $R_{21}$ is constructed in Fig. \ref{rys21-12}. It is of type  
$(N; N'; (N \rightarrow N); (N \rightarrow N); N''; N'''; N'''') \rightarrow Prop$, where socket $N$ corresponds to $n$, socket $N'$ to $o$, 
the first socket $(N \rightarrow N)$ to ${\tt node}^{in}$, 
the second socket $(N \rightarrow N)$ to ${\tt leaf}^{in}$, 
 socket $N''$ to $2$, socket $N'''$ to the first $1$, and socket $N''''$ to the second $1$.  
The relation $R_{21}(*; *; *; *; 2; 1; 1)$ of type $(N; N'; (N \rightarrow N); (N \rightarrow N)) \rightarrow Prop$, is the required condition $S_{21}$.  

~\\
For the second informal condition of $del$, and the corresponding operation  (denoted by $S_{22}$), the auxiliary operation $R_{22}$ is constructed in Fig. \ref{rys11-22}. It is of type $ ((N \rightarrow N); N'; N'';(N \rightarrow N); N''' ) \rightarrow Prop$, where the first socket $(N \rightarrow N)$ corresponds to ${\tt father}^{in}$, socket $N'$  to $o$,  socket $N''$ to $i$, the second socket $(N \rightarrow N)$ to ${\tt node}$, and socket $N'''$ to 1. 
So that $R_{22}({\tt father}^{in}; o; i; {\tt node}; 1))$ corresponds to  
$(Equal_N(o; i) + \neg Equal_N({\tt father}^{in}(i); {\tt father}^{in}(o)) +  \neg Equal_N(node^{in}(i); 1))$ denoted by $S({\tt father}; o; {\tt node}; i)$. 

By currying (applied to the socket $N''$ of $R_{22}$ corresponding to $i$) we get operation $R_{22}^c$ of type 
\\
$((N \rightarrow N); N'; (N \rightarrow N); N''' ) \rightarrow (N'' \rightarrow Prop)$. 

Now we are going to use operation $L^{\times}:((N''\rightarrow Prop); N'''')\rightarrow (N'''''\rightarrow Prop)$, see the previous section. 
Note that the socket $(N''\rightarrow Prop)$ of $L^{\times}$ is of the same type as the plug of $R_{22}^c$. So that apply amorphous composition to these two operations by making the connection between plug and the socket. The resulting operation denoted by $O$ is of type 
$( (N \rightarrow N); N'; (N \rightarrow N); N'''; N'''')\rightarrow (N'''''\rightarrow Prop)$
such that $O({\tt father}; o; {\tt node}; 1; n)(1)$ corresponds to 
  $(S({\tt father}; o; {\tt node}; 1) \ \ \times  \  \  S({\tt father}; o; {\tt node}; 2) \ \ \times ... \times \ \  S({\tt father}; o; {\tt node}; n))$. 
  
By applying uncurrying to the sub-socket $N'''''$ of the plug of operation $O$, we get operation $O^c$ of type $( (N \rightarrow N); N'; (N \rightarrow N); N'''; N''''; N''''')\rightarrow  Prop$, such that 
$O^c(*; *; *; 1; *;1)$ is the required condition $S_{22}$.

\begin{figure}[h]
	\centering
	\includegraphics[width=0.25\textwidth]{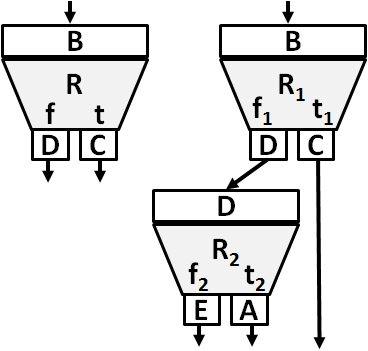}
	\caption{Operation ${\tt if\_then\_else}_{(B,C,D)}(R; t; f)$, and the  amorphous composition of ${\tt if\_then\_else}_{(B,C,D)}(R_1; t_1; f_1)$ and ${\tt if\_then\_else}_{(D,E,A)}(R_2; t_2; f_2)$ }
	\label{rys21}
\end{figure} 
%

%%%%%%%%
\subsection{Primitive operation ${\tt if\_then\_else}$ }
\label{if-then-else}

The {\bf if-then-else}, used to describe the operations $add$ and $del$, can not be constructed using the primitives introduced by now. It needs two operations $t:B \rightarrow C $ and $f:B \rightarrow D$, and a relation $R:B \rightarrow Prop$. For any $b:B$, it returns either $t(b)$ (if $R(b)$ is true) or $f(b)$ otherwise (else). 

Let the primitive operation  be denoted by 
${\tt if\_then\_else}_{(B,C,D)}$. It is of type $((B \rightarrow Prop); (B \rightarrow C); (B \rightarrow D)) \rightarrow (B  \rightarrow  (C || D))$. 
The type constructor $||$, i.e. exclusive disjunction, was introduced in Section  \ref{type-constructors}.

Application of the new primitive operation ${\tt if\_then\_else}_{(B,C,D)}$  to relation $R$, and two operations $t$ and $f$, i.e. ${\tt if\_then\_else}_{(B,C,D)}(R; t; f)$ results in the operation of type $B  \rightarrow  (C || D)$, see Fig. \ref{rys21}. 
That is, for any $b:B$, 
${\tt if\_then\_else}_{(B,C,D)}(R; t; f)(b)$  is either $t(b)$ if $R(b')$ is true, or $f(b)$ otherwise. Object $b'$ is a copy of $b$. 

If $f$ (or $t$) is $id_B$, then it means {\em do nothing}, i.e. return the input as the output. 

Amorphous composition of operation ${\tt if\_then\_else}_{(B,C,D)}(R_1; t_1; f_1)$ of type $B  \rightarrow  (C || D)$ and  operation ${\tt if\_then\_else}_{(D,E,A)}(R_2; t_2; f_2)$ of type $D  \rightarrow  (E || A)$ is shown in Fig. \ref{rys21}. The resulting operation is of type  $B  \rightarrow  (C || (E || A))$. 

A hardware interpretation of the new primitive operation requires a generic operation to evaluate conditions in their disjunctive normal form.

%%%%%%%
\subsection{ Constructions of $add$ and $del$}

Recall that $B$ denotes $(N;  N'; A; A'; A'')$ where $A$ denotes $(N \rightarrow N)$. Type $N$ corresponds the construction parameter $n$. Type $N'$ corresponds either to the father of the new node to be added, or to the node to be deleted. Type $A$ corresponds to operation ${\tt node}$. Type $A'$ corresponds to operation ${\tt father}$. Type $A''$ corresponds to operation ${\tt leaf}$. 

The operations $add$ and $del$ are supposed to be of type $B \rightarrow  ((B||B)||B)$.  Since the relations for the conditions of $add$ and $del$ have been already constructed, only the corresponding operations $t$ and $f$ are to be constructed.  
By introducing dumb sockets, corresponding to irrelevant parameters,  the relations become of type $B  \rightarrow Prop$, whereas the operations $t$ and $f$ are of type $B  \rightarrow B$.

\begin{figure}[h]
	\centering
	\includegraphics[width=0.98\textwidth]{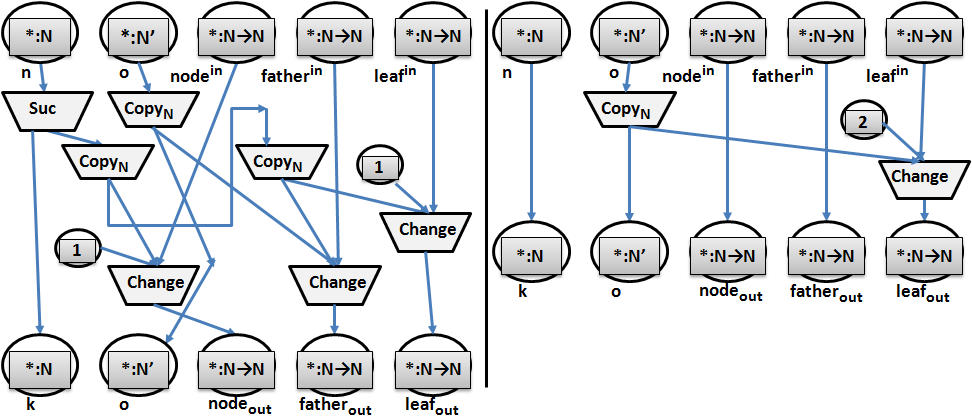}
	\caption{On the left, operation $f_{11}$ used to construction of $add$. On the right, operation $t_{12}$ used to construction of $add$ }
	\label{rys22a-22b}
\end{figure} 
\begin{figure}[h]
	\centering
	\includegraphics[width=0.98\textwidth]{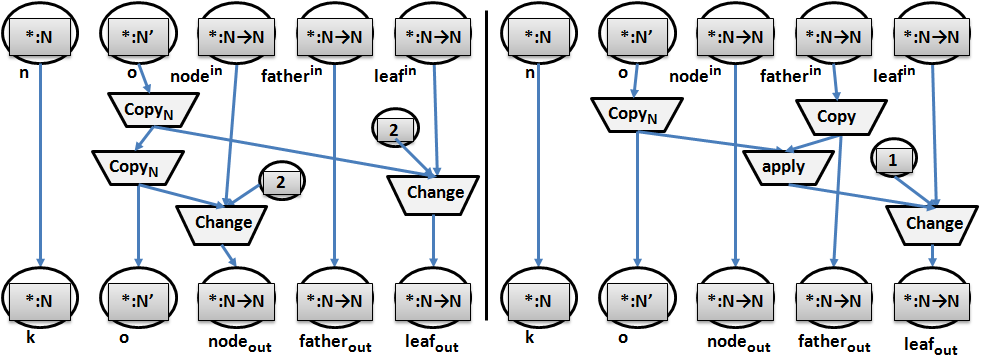}
	\caption{On the left, operation $f_{21}$  used to construction of $del$. On the right, operation  $t_{22}$ used to construction of $del$ }
	\label{rys23a-23b}
\end{figure} 

The operation $add$ is constructed as the amorphous composition of operation ${\bf if\_then\_else}_{(B,B,B)}(S_{11}; id_B; f_{11})$ of type $B \rightarrow (B'||B'')$,  and operation ${\bf if\_then\_else}_{(B,B,B)}(S_{12}; t_{12}; id_B)$  of type $B'' \rightarrow (B'''||B'''')$. The composition is done by connecting the plug $B''$ (corresponding to $f_{11}$) of the first operation to the socket of the second operation. The relations $S_{11}$ and $S_{12}$ have been already constructed whereas the operations $f_{11}$ and $t_{12}$ are constructed in Fig. \ref{rys22a-22b}. 

The operation $del$ is constructed as the amorphous composition of operation ${\bf if\_then\_else}_{(B,B,B)}(S_{21}; id_B; f_{21})$ of type $B \rightarrow (B'||B'')$,  and operation ${\bf if\_then\_else}_{(B,B,B)}(S_{22}; t_{22}; id_B)$ of type $B'' \rightarrow (B'''||B'''')$. The composition is done by connecting the plug $B''$ (corresponding to $f_{21}$) of the first operation to the socket of the second operation. The relations $S_{21}$ and $S_{22}$ have been already constructed  whereas the operations $f_{21}$ and $t_{22}$ are constructed in Fig. \ref{rys23a-23b}. 

% dodane 17 marca 2017 % 
Note that the operations $add$ and $del$ are of type $B \rightarrow  ((B||B)||B)$. Since the output consists of mutually disjoint plugs of the same type $B$, it seems that the plugs may be merged into one plug of type $B$, if it is needed. Then, the modified versions of $add$ and $del$ are  operations of type $B \rightarrow B$. However, then a special primitive type constructor should be introduced for merging type of the form $B||B$ into type $B$. 

Let us summarize the programming example.  
Once the operations $add$ and $del$ are constructed, they may be used to modify a tree data structure consisting of ${\tt node}$, ${\tt father}$, ${\tt leaf}$, and the parameter $n:N$ determining the current scope of the  structure.  This is done just by application of the operations $add$ and $del$ to an already existing tree structure. Perhaps it may be seen as an example of the pure function-level programming postulated by John Backus in his famous {\em 1977 ACM Turing Award Lecture} \cite{Backus}. 

Is this style of programming easy?  Since it is different from imperative and functional (based on term rewriting) programming, it may be quite hard for classical programmers. 

Note that the main part of the proposed programming style is in construction of functionals and relations. Once it has been done, the programming is easy. A program consists of applications and compositions of the functionals and the relations. The hardware interpretation suggests the compilation of the programs in tables of elementary operations and relations with dynamic link configurations between them. It seems that this style of programming is radically different from one corresponding to the von Neumann computer architecture, and from the one corresponding to term rewriting and lazy evaluation. 
%

%%%%%%%
%%%%%%%%%%%%%% 
\section{ Conclusion } 
% to też w skróconej wersji
The leitmotif of the paper was hardware interpretation (concrete grounding) of the introduced types, type constructors, and primitive generic operations that can be used to construct higher order objects (functionals and relations) in computing and programming. The constructions were merely described in an informal way. This may be an inspiration for creation a formal and technical specification  that may be useful for the design of electronic circuits.  
Actually, the idea is extremely simple, and consists in the management of links between plugs and sockets of elementary first order functions collected in huge arrays. 
Higher order functionals are efficient generic mechanisms for the management. 

It seems that the proposed approach may be considered as the second order intuitionistic Arithmetics. Along with the companion paper {\em Continuum as a primitive type} \cite{C} (that follows the idea of Brouwer) the work may shed some light on the grounding of the notion of functionals that is paramount for the Foundations of Mathematics. 
Finally let me express my personal view: Professor Luitzen E. J. Brouwer  was right and his intuitionism, when properly grasped and elaborated, constitutes the ultimate Foundation of Mathematics. % According to Brouwer  \cite{Brouwer1913}: ``''

% Bibliography
\bibliographystyle{ACM-Reference-Format}
\bibliography{biblio-TO-ANG}


\begin{thebibliography}{00}

%%% the Web address.

\ifx \showCODEN    \undefined \def \showCODEN     #1{\unskip}     \fi
\ifx \showDOI      \undefined \def \showDOI       #1{{\tt DOI:}\penalty0{#1}\ }
  \fi
\ifx \showISBNx    \undefined \def \showISBNx     #1{\unskip}     \fi
\ifx \showISBNxiii \undefined \def \showISBNxiii  #1{\unskip}     \fi
\ifx \showISSN     \undefined \def \showISSN      #1{\unskip}     \fi
\ifx \showLCCN     \undefined \def \showLCCN      #1{\unskip}     \fi
\ifx \shownote     \undefined \def \shownote      #1{#1}          \fi
\ifx \showarticletitle \undefined \def \showarticletitle #1{#1}   \fi
\ifx \showURL      \undefined \def \showURL       #1{#1}          \fi
% invisible to TeX
\providecommand\bibfield[2]{#2}
\providecommand\bibinfo[2]{#2}
\providecommand\natexlab[1]{#1}
\providecommand\showeprint[2][]{arXiv:#2}

\bibitem[\protect\citeauthoryear{Adamek, Milius, and Velebil}{Adamek
  et~al\mbox{.}}{2011}]%
        {Adamek}
\bibfield{author}{\bibinfo{person}{J. Adamek}, \bibinfo{person}{S. Milius},
  {and} \bibinfo{person}{J. Velebil}.} \bibinfo{year}{2011}\natexlab{}.
\newblock \showarticletitle{Semantics of higher-order recursion schemes}.
\newblock \bibinfo{journal}{{\em Logical Methods in Computer Science\/}}
  \bibinfo{volume}{7}, \bibinfo{number}{1:15} (\bibinfo{year}{2011}),
  \bibinfo{pages}{1--43}.
\newblock
\showDOI{%
\url{http://dx.doi.org/10.2168/LMCS-7 (1:15) 2011}}


\bibitem[\protect\citeauthoryear{Ambroszkiewicz}{Ambroszkiewicz}{2015}]%
        {C}
\bibfield{author}{\bibinfo{person}{Stanislaw Ambroszkiewicz}.}
  \bibinfo{year}{2015}\natexlab{}.
\newblock \bibinfo{title}{Continuum as a primitive type}.
\newblock \bibinfo{howpublished}{arxiv.org/abs/1510.02787}.
  (\bibinfo{year}{2015}).
\newblock
\showURL{%
\url{arxiv.org/abs/1510.02787}}


\bibitem[\protect\citeauthoryear{Backus}{Backus}{1978}]%
        {Backus}
\bibfield{author}{\bibinfo{person}{John Backus}.}
  \bibinfo{year}{1978}\natexlab{}.
\newblock \showarticletitle{Can Programming Be Liberated from the Von Neumann
  Style?: A Functional Style and Its Algebra of Programs}.
\newblock \bibinfo{journal}{{\em Commun. ACM\/}} \bibinfo{volume}{21},
  \bibinfo{number}{8} (\bibinfo{date}{Aug.} \bibinfo{year}{1978}),
  \bibinfo{pages}{613--641}.
\newblock
\showISSN{0001-0782}
\showDOI{%
\url{http://dx.doi.org/10.1145/359576.359579}}


\bibitem[\protect\citeauthoryear{Banach and Mazur}{Banach and Mazur}{1937}]%
        {BanachMazur1937}
\bibfield{author}{\bibinfo{person}{S. Banach} {and} \bibinfo{person}{S.
  Mazur}.} \bibinfo{year}{1937}\natexlab{}.
\newblock \showarticletitle{Sur les fonctions calculables}.
\newblock \bibinfo{journal}{{\em Annales de la Societe Polonaise de
  Mathematique\/}}  \bibinfo{volume}{16} (\bibinfo{year}{1937}),
  \bibinfo{pages}{223}.
\newblock


\bibitem[\protect\citeauthoryear{Bove and Capretta}{Bove and Capretta}{2005}]%
        {Bove}
\bibfield{author}{\bibinfo{person}{Ana Bove} {and} \bibinfo{person}{Venanzio
  Capretta}.} \bibinfo{year}{2005}\natexlab{}.
\newblock \showarticletitle{Modelling general recursion in type theory}.
\newblock \bibinfo{journal}{{\em Math. Struct. in Comp. Science\/}}
  \bibinfo{volume}{15} (\bibinfo{year}{2005}), \bibinfo{pages}{671--708}.
\newblock
\showDOI{%
\url{http://dx.doi.org/10.1017/S0960129505004822}}


\bibitem[\protect\citeauthoryear{Capalija and Abdelrahman}{Capalija and
  Abdelrahman}{2014}]%
        {capalija2014tile}
\bibfield{author}{\bibinfo{person}{Davor Capalija} {and}
  \bibinfo{person}{Tarek~S Abdelrahman}.} \bibinfo{year}{2014}\natexlab{}.
\newblock \showarticletitle{Tile-based bottom-up compilation of custom
  mesh-of-functional-units FPGA overlays}. In \bibinfo{booktitle}{{\em 2014
  24th International Conference on Field Programmable Logic and Applications
  (FPL)}}. IEEE, \bibinfo{pages}{1--8}.
\newblock


\bibitem[\protect\citeauthoryear{Cong, Huang, Ma, Xiao, and Zhou}{Cong
  et~al\mbox{.}}{2014}]%
        {cong2014fully}
\bibfield{author}{\bibinfo{person}{Jason Cong}, \bibinfo{person}{Hui Huang},
  \bibinfo{person}{Chiyuan Ma}, \bibinfo{person}{Bingjun Xiao}, {and}
  \bibinfo{person}{Peipei Zhou}.} \bibinfo{year}{2014}\natexlab{}.
\newblock \showarticletitle{A fully pipelined and dynamically composable
  architecture of CGRA}. In \bibinfo{booktitle}{{\em Field-Programmable Custom
  Computing Machines (FCCM), 2014 IEEE 22nd Annual International Symposium
  on}}. IEEE, \bibinfo{pages}{9--16}.
\newblock


\bibitem[\protect\citeauthoryear{Coquand}{Coquand}{2014}]%
        {Coq}
\bibfield{author}{\bibinfo{person}{T. Coquand}.}
  \bibinfo{year}{2014}\natexlab{}.
\newblock \bibinfo{title}{Coq Proof Assistant. Chapter 4 Calculus of Inductive
  Constructions}.
\newblock \bibinfo{howpublished}{www
  \url{http://coq.inria.fr/doc/Reference-Manual006.html}}.
  (\bibinfo{year}{2014}).
\newblock
\showURL{%
\url{http://coq.inria.fr/doc/Reference-Manual006.html}}
\newblock
\shownote{Site on www.}


\bibitem[\protect\citeauthoryear{Coquand and Huet}{Coquand and Huet}{1986}]%
        {Coquand1986}
\bibfield{author}{\bibinfo{person}{T. Coquand} {and}
  \bibinfo{person}{G\'{e}rard Huet}.} \bibinfo{year}{1986}\natexlab{}.
\newblock \bibinfo{booktitle}{{\em {The calculus of constructions}}}.
\newblock \bibinfo{type}{{T}echnical {R}eport} RR-0530.
  \bibinfo{institution}{INRIA}.
\newblock
\showURL{%
\url{http://hal.inria.fr/inria-00076024}}


\bibitem[\protect\citeauthoryear{Coquand and Huet}{Coquand and Huet}{1988}]%
        {Coquand1988}
\bibfield{author}{\bibinfo{person}{Thierry Coquand} {and}
  \bibinfo{person}{Gerard Huet}.} \bibinfo{year}{1988}\natexlab{}.
\newblock \showarticletitle{The Calculus of Constructions}.
\newblock \bibinfo{journal}{{\em Inf. Comput.\/}} \bibinfo{volume}{76},
  \bibinfo{number}{2-3} (\bibinfo{date}{Feb.} \bibinfo{year}{1988}),
  \bibinfo{pages}{95--120}.
\newblock
\showISSN{0890-5401}
\showDOI{%
\url{http://dx.doi.org/10.1016/0890-5401(88)90005-3}}


\bibitem[\protect\citeauthoryear{Curry}{Curry}{1964}]%
        {Curry}
\bibfield{author}{\bibinfo{person}{Haskell~B. Curry}.}
  \bibinfo{year}{1964}\natexlab{}.
\newblock \showarticletitle{Combinatory recursive objects of all finite types}.
\newblock \bibinfo{journal}{{\em Bull. Amer. Math. Soc.\/}}
  \bibinfo{volume}{70}, \bibinfo{number}{6} (\bibinfo{year}{1964}),
  \bibinfo{pages}{814--817}.
\newblock
\newblock
\shownote{www \url{http://projecteuclid.org/euclid.bams/1183526340}.}


\bibitem[\protect\citeauthoryear{De~Sutter, Raghavan, and Lambrechts}{De~Sutter
  et~al\mbox{.}}{2013}]%
        {de2013coarse}
\bibfield{author}{\bibinfo{person}{Bjorn De~Sutter}, \bibinfo{person}{Praveen
  Raghavan}, {and} \bibinfo{person}{Andy Lambrechts}.}
  \bibinfo{year}{2013}\natexlab{}.
\newblock \showarticletitle{Coarse-grained reconfigurable array architectures}.
\newblock In \bibinfo{booktitle}{{\em Handbook of signal processing systems}}.
  \bibinfo{publisher}{Springer}, \bibinfo{pages}{553--592}.
\newblock


\bibitem[\protect\citeauthoryear{Ershov}{Ershov}{1972a}]%
        {Ershov1972}
\bibfield{author}{\bibinfo{person}{Yu.L. Ershov}.}
  \bibinfo{year}{1972}\natexlab{a}.
\newblock \showarticletitle{Computable functionals of finite type}.
\newblock \bibinfo{journal}{{\em Algebra i Logika, English translation in
  Algebra and Logic, vol. 11 (1972), 203--242, AMS.\/}} \bibinfo{volume}{11},
  \bibinfo{number}{4} (\bibinfo{year}{1972}), \bibinfo{pages}{203--277}.
\newblock


\bibitem[\protect\citeauthoryear{Ershov}{Ershov}{1972b}]%
        {Ershov1973b}
\bibfield{author}{\bibinfo{person}{Yu.L. Ershov}.}
  \bibinfo{year}{1972}\natexlab{b}.
\newblock \showarticletitle{The theory of A-spaces}.
\newblock \bibinfo{journal}{{\em Algebra i Logika, English translation in
  Algebra and Logic, vol. 12 (1973), 209--232, AMS.\/}} \bibinfo{volume}{12},
  \bibinfo{number}{4} (\bibinfo{year}{1972}), \bibinfo{pages}{369--416}.
\newblock


\bibitem[\protect\citeauthoryear{Gammie}{Gammie}{2013}]%
        {gammie}
\bibfield{author}{\bibinfo{person}{Peter Gammie}.}
  \bibinfo{year}{2013}\natexlab{}.
\newblock \showarticletitle{Synchronous digital circuits as functional
  programs}.
\newblock \bibinfo{journal}{{\em ACM Computing Surveys (CSUR)\/}}
  \bibinfo{volume}{46}, \bibinfo{number}{2} (\bibinfo{year}{2013}),
  \bibinfo{pages}{21}.
\newblock


\bibitem[\protect\citeauthoryear{Girard}{Girard}{1971}]%
        {Girard71}
\bibfield{author}{\bibinfo{person}{J.Y. Girard}.}
  \bibinfo{year}{1971}\natexlab{}.
\newblock \showarticletitle{Une extension de l'interpretation de Godel a
  l'analyse, et son application a l'elimination des coupures dans l'analyse et
  dans la theorie des types}. In \bibinfo{booktitle}{{\em Proceedings of the
  Second Scandinavian Logic Symposium 1971}},
  \bibfield{editor}{\bibinfo{person}{J.~E. Fenstad}} (Ed.).
  \bibinfo{publisher}{North-Holland, Amsterdam}, \bibinfo{pages}{63--92}.
\newblock


\bibitem[\protect\citeauthoryear{Girard, Lafont, and Taylor}{Girard
  et~al\mbox{.}}{1989}]%
        {Girard1989}
\bibfield{author}{\bibinfo{person}{Jean-Yves Girard}, \bibinfo{person}{Yves
  Lafont}, {and} \bibinfo{person}{Paul Taylor}.}
  \bibinfo{year}{1989}\natexlab{}.
\newblock \bibinfo{booktitle}{{\em Proofs and Types}}.
  Vol.~\bibinfo{volume}{7}.
\newblock \bibinfo{publisher}{Cambridge University Press (Cambridge Tracts in
  Theoretical Computer Science)}.
\newblock
\showISBNx{0 521 37181 3}


\bibitem[\protect\citeauthoryear{G\"{o}del}{G\"{o}del}{1958}]%
        {Godel-Dialectica}
\bibfield{author}{\bibinfo{person}{K. G\"{o}del}.}
  \bibinfo{year}{1958}\natexlab{}.
\newblock \showarticletitle{{\"{U}}ber eine bisher noch nicht ben\"utzte
  Erweiterung des finiten Standpunktes}.
\newblock \bibinfo{journal}{{\em Dialectica\/}}  \bibinfo{volume}{10}
  (\bibinfo{year}{1958}), \bibinfo{pages}{280 -- 287}.
\newblock


\bibitem[\protect\citeauthoryear{Grzegorczyk}{Grzegorczyk}{1955a}]%
        {Grzegorczyk1955a}
\bibfield{author}{\bibinfo{person}{A. Grzegorczyk}.}
  \bibinfo{year}{1955}\natexlab{a}.
\newblock \showarticletitle{Computable functionals}.
\newblock \bibinfo{journal}{{\em Fundamenta Mathematicae\/}}
  \bibinfo{volume}{42} (\bibinfo{year}{1955}), \bibinfo{pages}{168--202}.
\newblock


\bibitem[\protect\citeauthoryear{Grzegorczyk}{Grzegorczyk}{1955b}]%
        {Grzegorczyk1955b}
\bibfield{author}{\bibinfo{person}{A. Grzegorczyk}.}
  \bibinfo{year}{1955}\natexlab{b}.
\newblock \showarticletitle{On the definition of computable functionals}.
\newblock \bibinfo{journal}{{\em Fundamenta Mathematicae\/}}
  \bibinfo{volume}{42} (\bibinfo{year}{1955}), \bibinfo{pages}{232--239}.
\newblock


\bibitem[\protect\citeauthoryear{Grzegorczyk}{Grzegorczyk}{1957}]%
        {Grzegorczyk1957}
\bibfield{author}{\bibinfo{person}{A. Grzegorczyk}.}
  \bibinfo{year}{1957}\natexlab{}.
\newblock \showarticletitle{On the definitions of computable real continuous
  functions}.
\newblock \bibinfo{journal}{{\em Fundamenta Mathematicae\/}}
  \bibinfo{volume}{44} (\bibinfo{year}{1957}), \bibinfo{pages}{61--71}.
\newblock


\bibitem[\protect\citeauthoryear{Grzegorczyk}{Grzegorczyk}{1964}]%
        {Grzegorczyk1964}
\bibfield{author}{\bibinfo{person}{A. Grzegorczyk}.}
  \bibinfo{year}{1964}\natexlab{}.
\newblock \showarticletitle{Recursive objects in all finite types}.
\newblock \bibinfo{journal}{{\em Fundamenta Mathematicae\/}}
  \bibinfo{volume}{54} (\bibinfo{year}{1964}), \bibinfo{pages}{73--93}.
\newblock


\bibitem[\protect\citeauthoryear{Hofmann}{Hofmann}{1999}]%
        {hofmann1999semantical}
\bibfield{author}{\bibinfo{person}{Martin Hofmann}.}
  \bibinfo{year}{1999}\natexlab{}.
\newblock \showarticletitle{Semantical analysis of higher-order abstract
  syntax}. In \bibinfo{booktitle}{{\em Proceedings. 14th Symposium on Logic in
  Computer Science}}. IEEE Computer Society, \bibinfo{pages}{204--204}.
\newblock


\bibitem[\protect\citeauthoryear{Jain, Li, Fahmy, and Maskell}{Jain
  et~al\mbox{.}}{2016}]%
        {jain2016adapting}
\bibfield{author}{\bibinfo{person}{Abhishek~Kumar Jain},
  \bibinfo{person}{Xiangwei Li}, \bibinfo{person}{Suhaib~A Fahmy}, {and}
  \bibinfo{person}{Douglas~L Maskell}.} \bibinfo{year}{2016}\natexlab{}.
\newblock \showarticletitle{Adapting the DySER architecture with DSP blocks as
  an Overlay for the Xilinx Zynq}.
\newblock \bibinfo{journal}{{\em ACM SIGARCH Computer Architecture News\/}}
  \bibinfo{volume}{43}, \bibinfo{number}{4} (\bibinfo{year}{2016}),
  \bibinfo{pages}{28--33}.
\newblock


\bibitem[\protect\citeauthoryear{Kleene}{Kleene}{1959a}]%
        {Kleene1959a}
\bibfield{author}{\bibinfo{person}{S.~C. Kleene}.}
  \bibinfo{year}{1959}\natexlab{a}.
\newblock \showarticletitle{Countable functionals}.
\newblock \bibinfo{journal}{{\em Constructivity in Mathematics: Proceedings of
  the colloquium held at Amsterdam\/}} (\bibinfo{year}{1959}),
  \bibinfo{pages}{81--100}.
\newblock


\bibitem[\protect\citeauthoryear{Kleene}{Kleene}{1959b}]%
        {Kleene1959b}
\bibfield{author}{\bibinfo{person}{S.~C. Kleene}.}
  \bibinfo{year}{1959}\natexlab{b}.
\newblock \showarticletitle{Recursive functionals and quantifiers of finite
  types I}.
\newblock \bibinfo{journal}{{\it Trans. Amer. Math. Soc.}}
  \bibinfo{volume}{91} (\bibinfo{year}{1959}), \bibinfo{pages}{1--52}.
\newblock


\bibitem[\protect\citeauthoryear{Kleene}{Kleene}{1963}]%
        {Kleene1963}
\bibfield{author}{\bibinfo{person}{S.~C. Kleene}.}
  \bibinfo{year}{1963}\natexlab{}.
\newblock \showarticletitle{Recursive functionals and quantifiers of finite
  types II}.
\newblock \bibinfo{journal}{{\it Trans. Amer. Math. Soc.}}
  \bibinfo{volume}{108} (\bibinfo{year}{1963}), \bibinfo{pages}{106--142}.
\newblock


\bibitem[\protect\citeauthoryear{Kreisel}{Kreisel}{1959}]%
        {Kreisel1959}
\bibfield{author}{\bibinfo{person}{G. Kreisel}.}
  \bibinfo{year}{1959}\natexlab{}.
\newblock \showarticletitle{Interpretation of analysis by means of functionals
  of finite type,}.
\newblock \bibinfo{journal}{{\em Constructivity in Mathematics: Proceedings of
  the colloquium held at Amsterdam\/}} (\bibinfo{year}{1959}),
  \bibinfo{pages}{101--128}.
\newblock


\bibitem[\protect\citeauthoryear{Lacombe}{Lacombe}{1955}]%
        {Lacombe1955}
\bibfield{author}{\bibinfo{person}{D. Lacombe}.}
  \bibinfo{year}{1955}\natexlab{}.
\newblock \showarticletitle{Remarques sur les operateurs recursifs et sur les
  fonctions recursives d'une variable reelle}.
\newblock \bibinfo{journal}{{\em Comptes Rendus de l'Academie des Sciences,
  Paris\/}}  \bibinfo{volume}{241} (\bibinfo{year}{1955}),
  \bibinfo{pages}{1250--1252}.
\newblock


\bibitem[\protect\citeauthoryear{Longley}{Longley}{2005}]%
        {Longley05notionsof}
\bibfield{author}{\bibinfo{person}{John Longley}.}
  \bibinfo{year}{2005}\natexlab{}.
\newblock \showarticletitle{Notions of computability at higher types I}. In
  \bibinfo{booktitle}{{\em In Logic Colloquium 2000}},
  \bibfield{editor}{\bibinfo{person}{R.~Cori}, \bibinfo{person}{A.~Razborov},
  \bibinfo{person}{S.~Todorcevic}, {and} \bibinfo{person}{C.~Wood}} (Eds.).
  \bibinfo{publisher}{Lecture Notes in Logic 19, ASL},
  \bibinfo{pages}{32--142}.
\newblock


\bibitem[\protect\citeauthoryear{Longley and Normann}{Longley and
  Normann}{2015}]%
        {longley2015higher}
\bibfield{author}{\bibinfo{person}{John Longley} {and} \bibinfo{person}{Dag
  Normann}.} \bibinfo{year}{2015}\natexlab{}.
\newblock \bibinfo{booktitle}{{\em Higher-Order Computability}}.
\newblock \bibinfo{publisher}{Springer}.
\newblock
\showISBNx{978-3-662-47991-9}
\showDOI{%
\url{http://dx.doi.org/10.1007/978-3-662-47992-6}}


\bibitem[\protect\citeauthoryear{Lyke, Christodoulou, Vera, and Edwards}{Lyke
  et~al\mbox{.}}{2015}]%
        {lyke2015introduction}
\bibfield{author}{\bibinfo{person}{James~C Lyke}, \bibinfo{person}{Christos~G
  Christodoulou}, \bibinfo{person}{G~Alonzo Vera}, {and}
  \bibinfo{person}{Arthur~H Edwards}.} \bibinfo{year}{2015}\natexlab{}.
\newblock \showarticletitle{An Introduction to Reconfigurable Systems}.
\newblock \bibinfo{journal}{{\it Proc. IEEE}} \bibinfo{volume}{103},
  \bibinfo{number}{3} (\bibinfo{year}{2015}), \bibinfo{pages}{291--317}.
\newblock


\bibitem[\protect\citeauthoryear{Ma, Aklah, and Andrews}{Ma
  et~al\mbox{.}}{2016}]%
        {ma2016just}
\bibfield{author}{\bibinfo{person}{Sen Ma}, \bibinfo{person}{Zeyad Aklah},
  {and} \bibinfo{person}{David Andrews}.} \bibinfo{year}{2016}\natexlab{}.
\newblock \showarticletitle{Just In Time Assembly of Accelerators}. In
  \bibinfo{booktitle}{{\em Proceedings of the 2016 ACM/SIGDA International
  Symposium on Field-Programmable Gate Arrays}}. ACM,
  \bibinfo{pages}{173--178}.
\newblock


\bibitem[\protect\citeauthoryear{Martin-L{\"{o}}f}{Martin-L{\"{o}}f}{1973}]%
        {Martin-Lof}
\bibfield{author}{\bibinfo{person}{P. Martin-L{\"{o}}f}.}
  \bibinfo{year}{1973}\natexlab{}.
\newblock \showarticletitle{An intuitionistic theory of types: predicative
  part}. In \bibinfo{booktitle}{{\em Logic Colloqium 1973}},
  \bibfield{editor}{\bibinfo{person}{H.~E. Rose} {and} \bibinfo{person}{J.~C.
  Shepherdson}} (Eds.). \bibinfo{publisher}{North-Holland, Amsterdam}.
\newblock


\bibitem[\protect\citeauthoryear{Mazur}{Mazur}{1963}]%
        {Mazur1963}
\bibfield{author}{\bibinfo{person}{S. Mazur}.} \bibinfo{year}{1963}\natexlab{}.
\newblock \showarticletitle{Computable analysis}.
\newblock \bibinfo{journal}{{\em Rozprawy Matematyczne\/}}
  \bibinfo{volume}{33} (\bibinfo{year}{1963}).
\newblock


\bibitem[\protect\citeauthoryear{Milner}{Milner}{1977}]%
        {milner1977fully}
\bibfield{author}{\bibinfo{person}{Robin Milner}.}
  \bibinfo{year}{1977}\natexlab{}.
\newblock \showarticletitle{Fully abstract models of typed $\lambda$-calculi}.
\newblock \bibinfo{journal}{{\em Theoretical Computer Science\/}}
  \bibinfo{volume}{4}, \bibinfo{number}{1} (\bibinfo{year}{1977}),
  \bibinfo{pages}{1--22}.
\newblock


\bibitem[\protect\citeauthoryear{Niedermeier, Kuper, and Smit}{Niedermeier
  et~al\mbox{.}}{2014}]%
        {niedermeier2014dataflow}
\bibfield{author}{\bibinfo{person}{Anja Niedermeier}, \bibinfo{person}{Jan
  Kuper}, {and} \bibinfo{person}{Gerard~JM Smit}.}
  \bibinfo{year}{2014}\natexlab{}.
\newblock \showarticletitle{A dataflow inspired programming paradigm for
  coarse-grained reconfigurable arrays}. In \bibinfo{booktitle}{{\em
  International Symposium on Applied Reconfigurable Computing}}. Springer,
  \bibinfo{pages}{275--282}.
\newblock


\bibitem[\protect\citeauthoryear{Palumbo, Sau, Fanni, Meloni, and
  Raffo}{Palumbo et~al\mbox{.}}{2016}]%
        {palumbo2016dataflow}
\bibfield{author}{\bibinfo{person}{Francesca Palumbo}, \bibinfo{person}{Carlo
  Sau}, \bibinfo{person}{Tiziana Fanni}, \bibinfo{person}{Paolo Meloni}, {and}
  \bibinfo{person}{Luigi Raffo}.} \bibinfo{year}{2016}\natexlab{}.
\newblock \showarticletitle{Dataflow-Based Design of Coarse-Grained
  Reconfigurable Platforms}. In \bibinfo{booktitle}{{\em Signal Processing
  Systems (SiPS), 2016 IEEE International Workshop on}}. IEEE,
  \bibinfo{pages}{127--129}.
\newblock


\bibitem[\protect\citeauthoryear{Paulson}{Paulson}{1986}]%
        {Paulson}
\bibfield{author}{\bibinfo{person}{L.~C. Paulson}.}
  \bibinfo{year}{1986}\natexlab{}.
\newblock \showarticletitle{Constructing Recursion Operators in Intuitionistic
  Type Theory}.
\newblock \bibinfo{journal}{{\em J. Symbolic Computation\/}}
  \bibinfo{volume}{2} (\bibinfo{year}{1986}), \bibinfo{pages}{325--355}.
\newblock


\bibitem[\protect\citeauthoryear{P\'{e}ter}{P\'{e}ter}{1934}]%
        {Peter}
\bibfield{author}{\bibinfo{person}{R. P\'{e}ter}.}
  \bibinfo{year}{1934}\natexlab{}.
\newblock \showarticletitle{Uber den Zussammenhang der verschiedenen Begriffe
  der rekursiven Funktion}.
\newblock \bibinfo{journal}{{\it Math. Ann.}}  \bibinfo{volume}{110}
  (\bibinfo{year}{1934}), \bibinfo{pages}{612--632}.
\newblock


\bibitem[\protect\citeauthoryear{P\'{e}ter}{P\'{e}ter}{1951a}]%
        {Peter1}
\bibfield{author}{\bibinfo{person}{R. P\'{e}ter}.}
  \bibinfo{year}{1951}\natexlab{a}.
\newblock \showarticletitle{Probleme der Hilbertschen Theorie der hoheren
  Stufen von rekursiven Funktionen}.
\newblock \bibinfo{journal}{{\em Acta mathematica Academiae Scientarum
  Hungaricae\/}}  \bibinfo{volume}{2} (\bibinfo{year}{1951}),
  \bibinfo{pages}{247--274}.
\newblock


\bibitem[\protect\citeauthoryear{P\'{e}ter}{P\'{e}ter}{1951b}]%
        {Peter2}
\bibfield{author}{\bibinfo{person}{R. P\'{e}ter}.}
  \bibinfo{year}{1951}\natexlab{b}.
\newblock \bibinfo{booktitle}{{\em Rekursive Funktionen}}.
\newblock \bibinfo{publisher}{AkademischerVerlag, Budapest. English translation
  published as Recursive Functions, Academic Press, 1967}.
\newblock


\bibitem[\protect\citeauthoryear{Platek}{Platek}{1966}]%
        {Platek1966}
\bibfield{author}{\bibinfo{person}{R. Platek}.}
  \bibinfo{year}{1966}\natexlab{}.
\newblock \bibinfo{title}{Foundations of recursion theory, Ph.D. thesis,
  Stanford University}.
\newblock   (\bibinfo{year}{1966}).
\newblock


\bibitem[\protect\citeauthoryear{Plotkin}{Plotkin}{1977}]%
        {Plotkin1977}
\bibfield{author}{\bibinfo{person}{Gordon~D. Plotkin}.}
  \bibinfo{year}{1977}\natexlab{}.
\newblock \showarticletitle{LCF considered as a programming language}.
\newblock \bibinfo{journal}{{\em Theoretical computer science\/}}
  \bibinfo{volume}{5}, \bibinfo{number}{3} (\bibinfo{year}{1977}),
  \bibinfo{pages}{223--255}.
\newblock


\bibitem[\protect\citeauthoryear{Reynolds}{Reynolds}{1974}]%
        {Reynolds}
\bibfield{author}{\bibinfo{person}{J. Reynolds}.}
  \bibinfo{year}{1974}\natexlab{}.
\newblock \showarticletitle{Towards a Theory of Type Structure}. In
  \bibinfo{booktitle}{{\em Colloque sur la Programmation 1974}}. Paris, France,
  \bibinfo{pages}{408--425}.
\newblock


\bibitem[\protect\citeauthoryear{Sazonov}{Sazonov}{1976}]%
        {Sazonov1976}
\bibfield{author}{\bibinfo{person}{V~Yu Sazonov}.}
  \bibinfo{year}{1976}\natexlab{}.
\newblock \showarticletitle{Degrees of parallelism in computations}. In
  \bibinfo{booktitle}{{\em International Symposium on Mathematical Foundations
  of Computer Science}}. Springer, \bibinfo{pages}{517--523}.
\newblock


\bibitem[\protect\citeauthoryear{Sch\"{u}rmann, Despeyroux, and
  Pfenning}{Sch\"{u}rmann et~al\mbox{.}}{2001}]%
        {Schurmann}
\bibfield{author}{\bibinfo{person}{Carsten Sch\"{u}rmann},
  \bibinfo{person}{Jo\"{e}lle Despeyroux}, {and} \bibinfo{person}{Frank
  Pfenning}.} \bibinfo{year}{2001}\natexlab{}.
\newblock \showarticletitle{Fundamental Study. Primitive recursion for
  higher-order abstract syntax}.
\newblock \bibinfo{journal}{{\em Theoretical Computer Science\/}}
  \bibinfo{volume}{266} (\bibinfo{year}{2001}), \bibinfo{pages}{1--57}.
\newblock


\bibitem[\protect\citeauthoryear{Scott}{Scott}{1970}]%
        {Scott1970}
\bibfield{author}{\bibinfo{person}{Dana~S. Scott}.}
  \bibinfo{year}{1970}\natexlab{}.
\newblock \showarticletitle{Outline of a mathematical theory of computation}.
  In \bibinfo{booktitle}{{\em Technical Monograph PRG-2, Oxford University
  Computing Laboratory, Oxford, England}}.
\newblock


\bibitem[\protect\citeauthoryear{Scott}{Scott}{1993}]%
        {Scott1993}
\bibfield{author}{\bibinfo{person}{D.~S. Scott}.}
  \bibinfo{year}{1993}\natexlab{}.
\newblock \showarticletitle{{A type-theoretical alternative to ISWIM, CUCH,
  OWHY}}.
\newblock \bibinfo{journal}{{\em Theoretical Computer Science\/}}
  \bibinfo{volume}{121} (\bibinfo{year}{1993}), \bibinfo{pages}{411--440}.
\newblock
\newblock
\shownote{first written in 1969 and widely circulated in unpublished form since
  then.}


\bibitem[\protect\citeauthoryear{Sheeran}{Sheeran}{1984}]%
        {Sheeran1984}
\bibfield{author}{\bibinfo{person}{Mary Sheeran}.}
  \bibinfo{year}{1984}\natexlab{}.
\newblock \showarticletitle{muFP, a Language for VLSI Design}. In
  \bibinfo{booktitle}{{\em Proceedings of the 1984 ACM Symposium on LISP and
  Functional Programming}} {\em (\bibinfo{series}{LFP '84})}.
  \bibinfo{publisher}{ACM}, \bibinfo{address}{New York, NY, USA},
  \bibinfo{pages}{104--112}.
\newblock
\showISBNx{0-89791-142-3}
\showDOI{%
\url{http://dx.doi.org/10.1145/800055.802026}}


\bibitem[\protect\citeauthoryear{Soare}{Soare}{1999}]%
        {Soare99thehistory}
\bibfield{author}{\bibinfo{person}{Robert~I. Soare}.}
  \bibinfo{year}{1999}\natexlab{}.
\newblock \showarticletitle{The history and concept of computability}. In
  \bibinfo{booktitle}{{\em in Handbook of Computability Theory}}.
  \bibinfo{pages}{3--36}.
\newblock


\bibitem[\protect\citeauthoryear{Tessier, Pocek, and DeHon}{Tessier
  et~al\mbox{.}}{2015}]%
        {tessier2015reconfigurable}
\bibfield{author}{\bibinfo{person}{Russell Tessier}, \bibinfo{person}{Kenneth
  Pocek}, {and} \bibinfo{person}{Andre DeHon}.}
  \bibinfo{year}{2015}\natexlab{}.
\newblock \showarticletitle{Reconfigurable computing architectures}.
\newblock \bibinfo{journal}{{\it Proc. IEEE}} \bibinfo{volume}{103},
  \bibinfo{number}{3} (\bibinfo{year}{2015}), \bibinfo{pages}{332--354}.
\newblock


\bibitem[\protect\citeauthoryear{{Univalent Foundations Program}}{{Univalent
  Foundations Program}}{2013}]%
        {HoTT}
\bibfield{author}{\bibinfo{person}{The {Univalent Foundations Program}}.}
  \bibinfo{year}{2013}\natexlab{}.
\newblock \bibinfo{booktitle}{{\em Homotopy Type Theory: Univalent Foundations
  of Mathematics}}.
\newblock \bibinfo{publisher}{Institute for Advanced Study,
  \url{http://homotopytypetheory.org/book}}.
\newblock
\showURL{%
\url{http://homotopytypetheory.org/book}}


\bibitem[\protect\citeauthoryear{Voevodsky}{Voevodsky}{2014}]%
        {vv}
\bibfield{author}{\bibinfo{person}{Vladimir Voevodsky}.}
  \bibinfo{year}{2014}\natexlab{}.
\newblock \showarticletitle{The Origins and Motivations of Univalent
  Foundations}.
\newblock \bibinfo{journal}{{\em IAS - The Institute Letter. Summer 2014,
  Institute for Advanced Study, Princeton, NJ, USA\/}} (\bibinfo{year}{2014}),
  \bibinfo{pages}{8--9}.
\newblock


\end{thebibliography}

\end{document}